\documentclass[12pt]{amsart}
\usepackage{epsfig}
\newfont{\sheaf}{eusm10 scaled\magstep1}

\newcommand{\ra}{\ensuremath{\rightarrow}}

\def\eea{\end{eqnarray*}}
\def\bea{\begin{eqnarray*}}
\def\Bbb{\bf}
\def\RR{{\Bbb R}}

\def\sA{{\mathcal{A}}}
\def\sD{{\mathcal{D}}}
\def\sC{{\mathcal{C}}}

\def\de{{\delta}}
\def\De{{\Delta}}

\def\a{{\alpha}}
\def\be{{\beta}}
\def\Ga{{\Gamma}}

\def\ZZ{{\Bbb Z}}

\newcommand{\Proof}{{\it Proof. }}
\newcommand{\QED}{{\hfill $Q.E.D.$}}

\newtheorem{teo}{Theorem}[section]
\newtheorem{df}[teo]{Definition}
\newtheorem{lem}[teo]{Lemma}
\newtheorem{cor}[teo]{Corollary}

\newtheorem{oss}[teo]{Remark}
\newtheorem{prop}[teo]{Proposition}

\newcommand{\C}{\ensuremath{\mathbb{C}}}
\newcommand{\LL}{\ensuremath{\mathbb{L}}}
\newcommand{\R}{\ensuremath{\mathbb{R}}}
\newcommand{\Z}{\ensuremath{\mathbb{Z}}}

\newcommand{\N}{\ensuremath{\mathbb{N}}}
\newcommand{\hol}{\ensuremath{\mathcal{O}}}

\newcommand{\PP}{\ensuremath{\mathbb{P}}}

\newcommand{\BB}{\ensuremath{\mathcal{B}}}

\newcommand{\VV}{\ensuremath{\mathbb{V}}}
\newcommand{\MM}{\ensuremath{\mathfrak{M}}}
\newcommand{\NNN}{\ensuremath{\mathfrak{N}}}
\newcommand{\gS}{\ensuremath{\mathfrak{S}}}

\usepackage{color} 
\newcommand{\rot}[1]{{\color{red} #1}}
\newtheorem{lemma}[teo]{Lemma}
\newtheorem{thm}[teo]{Theorem}

\newcommand{\ON}{\operatorname}
\newcommand{\inv}{^{^{-1}}}
\newcommand{\invv}{^{-2}}
\newcommand{\br}{\ON{Br}}

\renewcommand{\a}{\alpha}
\renewcommand{\b}{\beta}
\newcommand{\cg}{\gamma}

\renewcommand{\d}{\delta}

\newcommand{\eps}{\epsilon}
\newcommand{\z}{\zeta}

\newcommand{\s}{\sigma}
\newcommand{\oo}{\omega}

\allowdisplaybreaks

\begin{document}

\title[18 \today\ Moduli spaces and braid monodromy types.]
{Moduli spaces and braid monodromy types of bidouble covers of the quadric.}

\author{Fabrizio Catanese - Michael L\"onne -Bronislaw Wajnryb\\
         }

\date{\today}
\maketitle

\begin{abstract}
Bidouble covers $\pi : S \ra Q: = \PP^1 \times \PP^1$ of the quadric are
parametrized by connected families depending on four positive integers
$a,b,c,d$. In the special case where $ b=d$ we call them
abc-surfaces.

Such a Galois covering $\pi$  admits a small  perturbation
yielding a general 4-tuple covering of $Q$ with branch curve $\De$,
and a natural Lefschetz fibration obtained from a small perturbation
of the composition $ p_1 \circ \pi$.

We prove a more general result implying
         that the braid monodromy factorization corresponding to $\De$
determines  the three integers $a,b,c$ in the case of abc-surfaces. We
introduce a new method in order to distinguish  factorizations
which are not stably equivalent.

This result is in sharp contrast with  a previous result
         of the first and third author,
showing that the mapping class group factorizations corresponding to
the respective natural Lefschetz pencils are equivalent for
abc-surfaces with the same values of $ a+c, b$.
This result hints at the possibility that abc-surfaces with fixed
values of $ a+c, b$,
although diffeomorphic but not deformation equivalent, might  be not
canonically symplectomorphic.

          \end{abstract}

%\vfill
%\pagebreak
\section*{Introduction}

Bidouble covers of the quadric are smooth projective complex surfaces $S$
endowed with a (finite) Galois
covering $ \pi : S \ra Q : = \PP^1 \times \PP^1$ with Galois group
$(\Z /2 \Z)^2$.

More concretely, they are defined by a single pair of equations

\begin{eqnarray*}
         & z^2= f_{(2a,2b)}(x_0 , x_1; y_0 , y_1)\\
         & w^2=g_{(2c,2d)}(x_0 , x_1; y_0 , y_1)\\
\end{eqnarray*}

where we shall assume for simplicity that $ a,b,c,d\in\N_{\geq3}$,
and the notation $f_{(2a,2b)}$
denotes that
$f$ is a bihomogeneous polynomial, homogeneous of degree $2a$ in the
variables $x$,
and of degree $2b$ in the variables $y$.

These surfaces are simply connected and minimal of general type,
and they  were introduced in \cite{cat84} in order to show that the
moduli spaces
$\MM_{\chi,K^2}$ of smooth minimal surfaces of general type $S$ with
$ \chi(S): = \chi (\hol_S) = \chi$, $K^2_S = K^2 $,
need not be equidimensional or irreducible.

Given in fact our four integers $ a,b,c,d\in\N_{\geq3}$,  considering
the so called natural deformations of these bidouble covers,
defined by equations \footnote{in the following formula,
a polynomial of negative degree is identically zero. }

\begin{eqnarray*}
         & z^2= f_{(2a,2b)}(x_0 , x_1; y_0 , y_1) + w \
\Phi _{(2a-c,2b-d)}(x_0 , x_1; y_0 , y_1)\\
         & w^2=g_{(2c,2d)}(x_0 , x_1; y_0 , y_1) + z \
\Psi _{(2c-a,2d-b)}(x_0 , x_1; y_0 , y_1)\\
\end{eqnarray*}
one defines a bigger open subset $\NNN'_{a,b,c,d}$ of the moduli space,
whose closure $\overline{\NNN'}_{a,b,c,d}$
is an irreducible component of $\MM_{\chi,K^2}$, where
         $ \chi = 1 + (a-1)(b-1) + (c-1)(d-1) + (a+c-1)(b+d-1)$, and $K^2 = 8
(a+c-2)(b+d-2)$.

In general our knowledge about the moduli spaces $\MM_{\chi,K^2}$ is
rather scanty,
even if we make the drastic restriction only to consider the subset
$\MM^{00}_{\chi,K^2}$ corresponding to isomorphism classes $[S]$ of
simply connected minimal surfaces. $\MM_{\chi,K^2}$ is a
quasi-projective variety,
with a finite number of connected components, and $\MM^{00}_{\chi,K^2}$
is a union of connected components of $\MM_{\chi,K^2}$.

For a surface $S$ with $ [S] \in \MM^{00}_{\chi,K^2}$ there are (according to
Freedman's theorem \cite{free}) at most two topological types,
according to the parity of the intersection form $q_S : H^2 (S, \Z) \ra \Z$.
If $q_S$ is even (i.e., $Im (q_S) \subset 2 \Z)$ then $S$ is orientedly
homeomorphic to
a connected sum of copies of $\PP^1 \times \PP^1$ and of a K3 surface $Y$
(possibly taken with reversed orientation); if $q_S$ is odd
($Im (q_S) \not\subset 2 \Z)$
then $S$ is orientedly homeomorphic to a connected sum of copies of $\PP^2$
and of  $\PP^2$ with reversed orientation.

Thus the topology of $S$ does not give sufficient information in order to
distinguish the connected components  $\NNN \subset \MM^{00}_{\chi,K^2}$.
Note that if $\NNN$ is such a connected component and $S, S'$
are minimal surfaces with $[S], [S'] \in \NNN$ (we then simply say that
$S, S'$  have the same deformation type), then, by the classical
theorem of Ehresmann (\cite{ehre}), there exists an orientation
preserving  diffeomorphism
$\Psi : S' \ra S$ carrying the canonical class $c_1 (K_S)$ to
the canonical class $c_1 (K_{S'}).$

In \cite{cat02} (see also \cite{cat09})it was shown that to each such surface
$S$ one can associate a  symplectic manifold $(S, \omega)$,
unique up to symplectomorphism, and such that the De Rham
class of $\omega$ equals $c_1 (K_S)$. Moreover, it was shown
that this symplectomorphism class,
called `canonical symplectomorphism class',  is an invariant of
the connected component $\NNN$.

The so called Manetti surfaces yield examples (see \cite{man4},
\cite{cat09}, \cite{catcime})
of surfaces $S, S'$ lying in distinct
connected components of $\MM_{\chi,K^2}$ but
having  the same canonical symplectomorphism class.

Whether  this also occurs for $\MM^{00}_{\chi,K^2}$ is an interesting
open question,
especially motivated by the results of \cite{c-w}, where the
         families  $\NNN_{a,b,c}$ of the so called `abc'-surfaces,
essentially a partial closure of $\NNN'_{a,b,c,b}$,
were shown to provide examples of simply connected surfaces $S, S'$ belonging
         to distinct
connected components ($\NNN_{a,b,c}$, resp. $\NNN_{a+k,b,c-k}$)
of $\MM^{00}_{\chi,K^2}$, yet diffeomorphic through a
         diffeomorphism  preserving the orientation and the canonical class.

The starting point was that a partial closure $\NNN_{a,b,c,d}$ of
$\NNN'_{a,b,c,d}$,
obtained by allowing the base to be a more general Segre-Hirzebruch surface,
and allowing the bidouble covers to have also Du Val
singularities, is, under suitable numerical conditions on
$a,b,c,d$, an irreducible connected component of  $\MM_{\chi,K^2}$.

Then it was proven in \cite{c-w} that if $S$ is an abc-surface,
and $S'$ is an $a'b'c'$-surface, then $S$ and $S'$ are diffeomorphic
if and only if $b = b', a+c = a' + c'$.

The diffeomorphism between $ S$ and $S'$ was obtained
as a consequence of a classical theorem of Kas
on the diffeomorphism type of differentiable  Lefschetz fibrations.

Here, the holomorphic map  $\tilde{\varphi} : S \ra \PP^1$, obtained as
the composition of $\pi : S \ra \PP^1 \times \PP^1$
with the first projection, admits a small differentiable
perturbation which yields a natural
symplectic Lefschetz fibration $\varphi : S \ra \PP^1$.

The bulk of the proof was to show that $\varphi$, $\varphi '$
are isomorphic differentiable Lefschetz fibrations according to the criterion
of Kas.

Kas \cite{kas} shows indeed that
the isomorphism class of a differentiable Lefschetz fibration
of genus $g\geq 2$ is completely determined by the equivalence class
of its monodromy factorization in the Mapping class group,
for the equivalence relation generated by Hurwitz equivalence
and by simultaneous conjugation.

Here, the Mapping class group monodromy factorization of a Lefschetz
fibration is the
sequence of  positive Dehn twists
associated to a quasi-basis of the fundamental group
$\pi_1(\PP^1-\{b_i\},b_0)$, whose conjugacy classes yield the local
monodromies.

The change of quasi-basis of the fundamental group
$\pi_1(\PP^1-\{b_i\},b_0)$ leads to the Hurwitz equivalence
of factorizations, while simultaneous conjugation accounts for
the possible different choices of
an (orientation preserving)
diffeomorphism of the fibre over the base point $b_0$
with a standard
Riemann surface of genus $g$.

On the other hand, the Mapping class group monodromy
factorization of this Lefschetz fibration is nothing else than the
homomorphic image
of the Braid Monodromy factorization
which corresponds to the
branch curve $\De \subset \PP^1 \times \PP^1$ of a symplectic
perturbation of $\pi$
and to its first projection onto $ \PP^1$, and which has factors in
the braid group of the sphere $\br_n$, where $n$ is the vertical
degree of $\De$:

$$
\br_n \quad = \quad
\left\langle \: \s_1,..., \s_{n-1}\:\left|
\begin{array}{l}
\s_i\s_j = \s_j\s_i,\quad \text{ if } |i-j|>1\\
\s_i\s_{i+1}\s_i = \s_{i+1}\s_i\s_{i+1}\\
\s_1\cdots\s_{n-1}\s_{n-1}\cdots\s_1=1
\end{array}
\right.\right\rangle
$$

Determining whether
these Braid Monodromy factorizations are equivalent,  or not,
is the clue to  deciding about the existence of
a diffeomorphism between $S, S'$  commuting with the Lefschetz fibrations
$\varphi, \varphi '$, and yielding a canonical
symplectomorphism.

Our main result says that such a diffeomorphism cannot exist,
by showing that the corresponding Braid Monodromy factorizations are
not m-equivalent, in the terminology introduced by
Auroux and Katzarkov. The relation of m-equivalence
was introduced by Auroux
and Katzarkov in \cite{a-k} in order to obtain invariants of
symplectic 4-manifolds. It is obtained
by allowing not only Hurwitz equivalence and simultaneous conjugation, but also
creation/cancellation of admissible
pairs of a positive and of a negative node,  where a node  (and
then also the  corresponding full twist $\beta$) is said
to be admissible if the inverse image of the node inside the
ramification divisor consists of two disjoint smooth branches.

We propose here to use the more informative name:
stable equivalence, instead of m-equivalence.

More precisely, we prove the following Main Theorem:

\begin {teo}
The braid monodromy factorizations associated to
an abc-surface $S$  and to
       an $a'b'c'$-surface  $S'$ are not stably equivalent, except in the
trivial cases
       $ a = a', b = b', c =  c'$ or $ a = c', b = b', c =  a'$ , or  \ $
a=c=b', a' = c' = b$.
\end{teo}

The result  relies on a complete description of the braid monodromy 
factorization class
associated to $\De$: this is given in Theorem \ref{BMF}, which proves
indeed much more than what we need
for the present purposes.

We state here a simpler byproduct of the cited result, namely:

\begin{thm}
\label{bmf}
There is a braid monodromy factorization of the curve $\De$
associated to a bidouble cover $S$ of type $a,b,c,d$
     whose braid monodromy group $H\subset\br_{4(b+d)}$ is generated,
unless we are in the cases  $$ {\rm(I) }\  c = 2a \ {\rm and } \  d =
2b \ , \ {\rm or\  (II)} \ a = 2c
\ {\rm and } \  b = 2d \,$$

   by
the following powers of half-twists:
$\s_{a_i},\s_{c_i}$ for $i=1,\dots 2b-1$, $\s_{b_\imath},\s_{d_\imath}$,
for $\imath=1,\dots,2d-1$,
$\s_{p_{2b}}^2,\s_{q_{2d}}^2,\s_{s},\s_{u'}^3,\s_{u''}^3$.

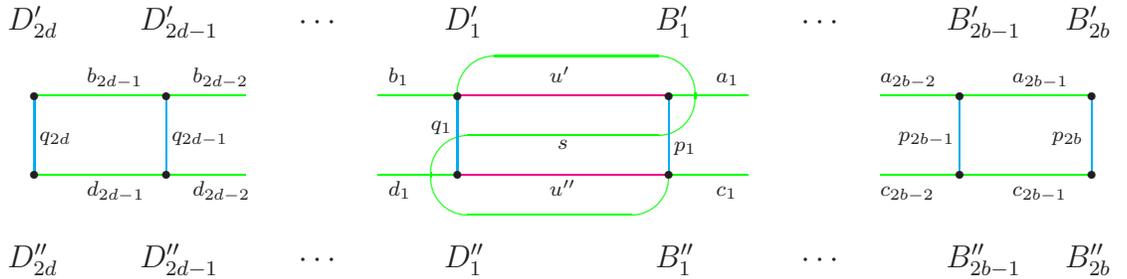
\begin{picture}(420,120)(-210,-60)

{\color{green}
\put(-200,-15){\line(1,0){80}}
\put(-70,-15){\line(1,0){30}}
\put(200,-15){\line(-1,0){80}}
\put(70,-15){\line(-1,0){30}}
\put(-200,15){\line(1,0){80}}
\put(-70,15){\line(1,0){30}}
\put(200,15){\line(-1,0){80}}
\put(70,15){\line(-1,0){30}}

\put(5,15){\oval(90,30)[t]}
\put(0,15){\oval(100,30)[rb]}
\put(0,-15){\oval(100,30)[lt]}
\put(-5,-15){\oval(90,30)[b]}
\color{magenta}
\put(-40,15){\line(1,0){80}}
\put(-40,-15){\line(1,0){80}}
\color{cyan}
\put(-200,-15){\line(0,1){30}}
\put(-150,-15){\line(0,1){30}}
\put(-40,-15){\line(0,1){30}}
\put(40,-15){\line(0,1){30}}
\put(150,-15){\line(0,1){30}}
\put(200,-15){\line(0,1){30}}
}

\put(40,15){\circle*{3}}
\put(150,15){\circle*{3}}
\put(200,15){\circle*{3}}
\put(-40,15){\circle*{3}}
\put(-150,15){\circle*{3}}
\put(-200,15){\circle*{3}}
\put(40,-15){\circle*{3}}
\put(150,-15){\circle*{3}}
\put(200,-15){\circle*{3}}
\put(-40,-15){\circle*{3}}
\put(-150,-15){\circle*{3}}
\put(-200,-15){\circle*{3}}

\put(-210,-50){$D''_{2d}$}
\put(-160,-50){$D''_{2d-1}$}
\put(-45,-50){$D''_{1}$}
\put(35,-50){$B''_{1}$}
\put(145,-50){$B''_{2b-1}$}
\put(190,-50){$B''_{2b}$}
\put(-210,40){$D'_{2d}$}
\put(-160,40){$D'_{2d-1}$}
\put(-45,40){$D'_{1}$}
\put(35,40){$B'_{1}$}
\put(145,40){$B'_{2b-1}$}
\put(190,40){$B'_{2b}$}

\put(-100,-50){$\cdots$}
\put(-100,40){$\cdots$}
\put(90,-50){$\cdots$}
\put(90,40){$\cdots$}

\put(-198,-2){$\scriptstyle q_{2d}$}
\put(-148,-2){$\scriptstyle q_{2d-1}$}
\put(-50,2){$\scriptstyle q_{1}$}
\put(42,-6){$\scriptstyle p_{1}$}
\put(127,-2){$\scriptstyle p_{2b-1}$}
\put(185,-2){$\scriptstyle p_{2b}$}

\put(-180,20){$\scriptstyle b_{2d-1}$}
\put(-140,20){$\scriptstyle b_{2d-2}$}
\put(-66,20){$\scriptstyle b_{1}$}
\put(-5,20){$\scriptstyle u'$}
\put(58,20){$\scriptstyle a_{1}$}
\put(120,20){$\scriptstyle a_{2b-2}$}
\put(170,20){$\scriptstyle a_{2b-1}$}
\put(-180,-23){$\scriptstyle d_{2d-1}$}
\put(-140,-23){$\scriptstyle d_{2d-2}$}
\put(-66,-23){$\scriptstyle d_{1}$}
\put(-5,-23){$\scriptstyle u''$}
\put(58,-23){$\scriptstyle c_{1}$}
\put(120,-23){$\scriptstyle c_{2b-2}$}
\put(170,-23){$\scriptstyle c_{2b-1}$}

\put(-2,-5){$\scriptstyle s$}
\end{picture}\\
The factorization is such that
\begin{enumerate}
\item
each ($\pm$)full-twist factor is of type $p$ or $q$,
\item
the weighted count of ($\pm$)full-twist factors of type $p$ yields
$8ab-2(ad+bc)$,
\item
the weighted count of ($\pm$)full-twist factors of type $q$ yields
$8cd-2(ad+bc)$.
\end{enumerate}
\end{thm}

The above result shows  that the braid monodromy group $H$ depends
only upon the numbers $ b $ and $d$, provided for instance that we have
nonvanishing  of the  respective numbers
$8ab-2(ad+bc)$,
$8cd-2(ad+bc)$,
or provided that we are in the case $b=d$.

A fortiori if the groups $H$ are the same for different choices of $(a,b,c,d)$
then the  fundamental groups $\pi_1 ( Q \setminus \De) $ are isomorphic.

This forces us to look more carefully into the problem of distinguishing
classes of factorizations, and for this reason, later on, we
introduce a  technical
novelty which consists  in finding a new effective method for disproving
stable-equivalence.

The method goes as follows: assume that we consider a group $G$
and a factorization of the identity in $G$
$$ \a_1 \circ   \dots  \circ \a_m = 1$$
and a set $\BB$ of elements $\be_j \in G$.
We  define stable-equivalence with respect to $\BB$ simply by
considering the equivalence relation generated by Hurwitz equivalence,
simultaneous conjugation and creation/cancellation of consecutive factors
$  \be_j \circ \be_j^{-1}$.

Then  {\bf refined  invariants } of the stable equivalence class of
the factorization are obtained
as follows:

1) let $H$ be the subgroup of $G$ generated by the $\a_i$'s
and let $\hat{H}$ be the subgroup of $G$ generated by the $\a_i$'s
and by the $\be_j$'s (in our case, $H$ will be called the monodromy group,
and $\hat{H}$ the stabilized monodromy group).

Then $\hat{H}$ is a first invariant.
\smallskip

2) Let $\sC$ be the set of conjugacy classes in the  group $H$,
and let $\hat{\sC}$ be the set of conjugacy classes in the  group
$\hat{H}$, so that we have a natural map $\sC \ra \hat{\sC}$,
such that $A \mapsto \hat{A}$.

Write $\hat{\sC}$ as a disjoint union $\hat{\sC}^+ \cup \hat{\sC}^-
\cup \hat{\sC}^0$,
where $\hat{\sC}^0$ is the set of conjugacy classes of elements $a$ which
are conjugate
to their inverse $ a^{-1}$, and where $\hat{\sC}^-$ is the set of the inverse
conjugacy classes of the classes in $\hat{\sC}^+$.

3) Associate to the factorization $ \a_1 \circ   \dots  \circ \a_m = 1$
the function $ s : \hat{\sC}^+ \ra \Z$ such that $ s(c)$, for $ c \in
\hat{\sC}^+$,
is the algebraic number of occurrences of $c$ in the sequence of
conjugacy classes
of the $ \a_j $'s (i.e., an occurrence of $c^{-1}$ counts as $-1$ for $ s(c)$).

The function $ s : \hat{\sC}^+ \ra \Z$ is  our second and most
important invariant.

\begin{oss}
The calculation of the function $ s : \hat{\sC}^+ \ra \Z$
      presupposes however a detailed knowledge of the group $  \hat{H}$.
For this reason it will be convenient to find some coarser
derived invariant.
\end{oss}

{\bf Substrategy I.}

Assume that  we can write $\{\a_1, \dots, \a_m\} $ as a disjoint
union $\sA_1 \cup \sA_2 \cup \sD \cup \sA'_1 \cup \sA'_2 $, such that
the set $\sA_j, \ j=1,2,$
is contained in a conjugacy class $A_j \subset H$,
the set $\sA'_j, \ j=1,2,$ is contained in the conjugacy class
$A_j^{-1} \subset H$.
Assume that the elements  in $\sA_1 \cup \sA_2$ are contained in a conjugacy
class $C$ in $G$ such that $ C \cap C^{-1} = \emptyset$ but that
the set $\sD$ is disjoint from the union of the two conjugacy classes of $G$,
$ C \cup C^{-1} $.

If we then prove that $\hat{A_1} \neq \hat{A_2} $ (this of course implies
$A_1 \neq A_2 $) we may assume w.l.o.g.  that
$\hat{A_1} , \hat{A_2} \in \hat{\sC}^+$, and  then the unordered pair
of positive numbers $( |s(\hat{A_1} )|, |s(\hat{A_2}| )$  is our derived
numerical invariant of the factorization (and can easily be calculated from the
cardinalities of the four sets  as $(||\sA_1| - |\sA'_1|| , ||\sA_2|
- |\sA'_2||)$.

\bigskip
{\bf Substrategy II.}

This is the strategy to show that $\hat{A_1} \neq \hat{A_2} $
and goes as follows.

Assume further that we have a subgroup $ \tilde{H} \supset \hat{H}$,
and a group homomorphism $\rho :  \tilde{H}  \ra \Sigma$ such that
$\rho ( \BB ) = 1$.  Assume also that the following key property holds.

{\bf Key property.}

      For each element $ \a_j \in A_j \subset H$ there exists
an element $ \tilde{\a}_j \in \tilde{H}$
      such that $$ \a_j = \tilde{\a}_j^2 ,$$
and moreover that this element is unique in $\tilde{H}$ (a fortiori, it
will suffice that this element is unique in $G$).

Proving that  $ \rho (\tilde{\a}_1)$ is not conjugate to
$ \rho (\tilde{\a}_2)$ under the action of $ \rho (H) = \rho (\hat{H})
\subset \Sigma$
shows finally that $\a_1$ is not conjugate to $\a_2$ in $\hat{H}$.

Since, if there is $ h \in \hat{H}$ such that $\a_1 = h^{-1} \a_2 h$,
then $\tilde{\a}_1^2  = \a_1 = (h^{-1} \tilde{\a}_2 h)^2$, whence
$\tilde{\a}_1  = h^{-1} \tilde{\a}_2 h$, a contradiction.

In our concrete case, we are able to determine the braid monodromy group $H$,
and we observe that, since we have a so called cuspidal factorization,
all the factors $\a_i$ belong to only four conjugacy classes
in the group $G$ (the classes of $\sigma_1$,$\sigma_1^2$,
$\sigma_1^{-2}$,$\sigma_1^3$
in the braid group), corresponding geometrically to vertical tangencies of the
branch curve, respectively nodes and cusps.
Three of these classes are positive and only one is negative (the one
of $\sigma_1^{-2}$).

We see that  the nodes belong to  conjugacy classes $A_1, A_2,
A_1^{-1}, A_2^{-1}
\subset H$, the positive nodes belonging to  $A_1 \cup  A_2$,
and we show, using a representation $\rho$ of  a certain subgroup $\tilde{H}$
of `liftable' braids (i.e., braids which centralize the monodromy
homomorphism, whence are liftable to the mapping class group of
the curve associated to the monodromy homomorphism) into a symplectic group
$\Sigma$ with $\Z / 2$ coefficients,
that the classes $\hat{A_1}, \hat{A_2}$ are distinct in $\hat{H}$. Moreover,
these are classes in  $\sC^+$, since these are positive classes in the
braid group.
We calculate then easily the above function for these two conjugacy classes,
i.e., the pair of numbers $ ( s(\hat{A_1}), s(\hat{A_2})) $.

This new method  and the above results represent the first
positive step towards the realization
of a more general program set up by Moishezon (\cite{mois1}, \cite{mois2})
in order to produce braid monodromy invariants which should distinguish
the connected components of a moduli space $\MM_{\chi,K^2}$.

Moishezon's program is based on the consideration (assume here for simplicity
that $K_S$ is ample) of
a general projection $\Psi_m : S \ra \PP^2$ of a pluricanonical
embedding $\Phi_m : S \ra \PP^{P_m -1}$, and of the braid monodromy
factorization
corresponding to the (cuspidal) branch curve $B_m$ of $\Psi_m$.

An invariant of the connected component is here given by the equivalence class
(for Hurwitz equivalence plus simultaneous conjugation) of this
braid monodromy factorization. Moishezon, and later Moishezon-Teicher,
calculated a coarser invariant, namely the fundamental group
$\pi_1(\PP^2- B_m)$.

This group turned out to be not overly complicated,
and in fact, as shown in many cases in \cite{a-d-k-y}, it tends  to give
no extra information beyond the one given by the topological
         invariants of $S$ (such as $\chi, K^2$).

Auroux and Katzarkov showed that, for $ m >>0$, the stable-equivalence class
of the above braid monodromy factorization determines the canonical
symplectomorphism class of $S$, and conversely.

As we have already remarked, in the case of abc-surfaces
% with $ c \neq 3a$, and $ a \neq 3c$,
the braid monodromy groups are determined
by $b$ (up to conjugation), hence the fundamental groups  $\pi_1 ( Q
\setminus \De) $ are isomorphic for a fixed value of $b$.

So, Moishezon's
technique produces no invariants.
%Moreover, even when $ c= 3a$ or $ a = 3c$
%the stabilized fundamental groups are the same, once $b$ is fixed.

Therefore our result indicates that one should try to go all the way
to understanding stable-equivalence classes of pluricanonical
braid monodromy factorizations.  Let us try  to describe here how
this could be done.

Define  $ p : S \ra \PP^2$ as the morphism  $p$
obtained as the composition of
$\pi : S \ra Q$ with the embedding $ Q \hookrightarrow \PP^3$
followed by a general projection  $\PP^3 \dashrightarrow \PP^2$.

In the even more special case of abc-surfaces such that $ a+c = 2b$,
         the m-th pluricanonical mapping $\Phi_m : S \ra \PP^{P_m -1}$
has a (non generic) projection given by the composition of $p$ with
a Fermat type map  $\nu_r : \PP^2  \ra \PP^2$
(given by $ \nu_r (x_0, x_1, x_2) = (x_0^r, x_1^r, x_2^r) $
in a suitable linear coordinate system),
where $ r : = m ( 2b -2) $.

Let $B$ be the branch curve of a generic perturbation of $ p$:
then the braid monodromy factorization
corresponding to $B$ can be calculated from the braid monodromy factorization
corresponding to $\De$.

We hope, in a sequel to this paper,
to be able to determine whether these braid monodromy factorizations are
equivalent, respectively stably-equivalent,  for abc-surfaces such that
$ a+c = 2b$.

The problem of calculating the braid monodromy factorization corresponding
to the (cuspidal) branch curve $B_m$ starting from the
braid monodromy factorization of $B$ has been addressed, in the
special case $ m=2$,
by Auroux and Katzarkov (\cite{a-k2}). Iteration of
their formulae should lead to the calculation of
         the braid monodromy factorization corresponding
to the (cuspidal) branch curve $B_m$ in the case, sufficient for applications,
         where $m$ is a sufficiently large power of $2$.

Here are the contents of the article.

The first section is devoted to making precise the concept of
a perturbation of $\pi$. In fact, it will be shown that a dianalytic
perturbation suffices, i.e., one corresponding to a covering given by
equations where
$\Phi, \Psi$ are either holomorphic or antiholomorphic in each
single variable $x$, $y$.

The invariants of the branch curve $\De$
(genus, degree, number of nodes, cusps and vertical tangents) are computed.

The second section is devoted to the description of
the braid monodromy group $H$ of $\De$,
and to the conjugacy classes in $H$ of the factors of the
braid monodromy factorization. The computation is based
on a degeneration of the branch curves to real curves $\{ f=0\},\{ g=0\} $
which are each the union
of horizontal  lines with the graph of a rational function: this type
of degeneration seems to be very suitable also for writing short proofs of
existing results.

The third section is based on the classical correspondence between
4-tuple covers and triple covers, given by the surjection
$\gS_4 \ra \gS_3$ whose kernel is the Klein group $(\Z / 2 )^2$.

We obtain a corresponding triple cover and a
resulting homomorphism of its  mapping class
group  to the symplectic group acting on the  $\Z / 2$  homology.

This homomorphism sends the 'extra' factors $\be_j$ to the identity
and transforms the Dehn twists $\alpha_i$ corresponding to the nodes
of $\De$ to Picard Lefschetz transformations  which are shown to
be not conjugated under the image of the monodromy group.

\section{Perturbed simple bidouble covers}

Consider the direct sum $\VV$ of two complex line bundles
$\LL_1 \oplus \LL_2$
on a compact complex manifold  $X$, and the subset  $Z$
of $\VV$ defined by the following pair of equations,

\begin{eqnarray*}
         & z^2= f(x) + w \
\Phi (x)\\
         & w^2=g(x) + z \
\Psi (x)\\
\end{eqnarray*}

where $f,g$ are respective holomorphic sections of the line bundles
$\LL_1^{\otimes 2},  \LL_2^{\otimes 2}$ (we shall also write
$ f \in H^0 (\hol_X(2 L_1)),  g \in H^0 (\hol_X(2 L_2))$, denoting by
$L_1, L_2$ the associated
Cartier divisors), and where, for the time being, $\Phi$ is a
differentiable section
of  $\LL_1^{\otimes 2} \otimes \LL_2^{\otimes -1}$,
$\Psi$ is a differentiable section
of  $\LL_2^{\otimes 2} \otimes \LL_1^{\otimes -1}$.

\begin{lem}\label{smallperturbation}
Assume that the two divisors $\{ f=0\}$ and $\{ g=0\}$ are smooth and
intersect transversally.
Then, for $|\Phi| < < 1$, $|\Psi| < < 1$, $Z$ is a smooth submanifold
of $\VV$, and the projection
$\pi : \VV \ra X$ induces  a finite covering  $ Z \ra X$ of degree
$4$, with ramification divisor (i.e., critical set) $ R : = \{ 4 zw =
\Phi \Psi \}$, and
with branch divisor (i.e., set of critical values) $\De = \{ \de(x) = 0\}$,
where $$ - \frac{1}{16} \de = - f^2 g^2 - \frac{9}{8} fg (\Phi
\Psi)^2 + (\Psi)^2 f^3 +
        (\Phi)^2 g^3 + \frac{27}{16^2} (\Phi \Psi)^4.$$
\end{lem}

\Proof
$Z$ is locally defined by two complex valued functions which are holomorphic
polynomials in $z$ and $w$. If we compactify the rank two vector bundle to a
$\PP^2_{\C}$-bundle, and homogenize the equations to
\begin{eqnarray*}
         & z^2= f(x) u^2 + w u \
\Phi (x)\\
         & w^2=g(x)u^2   + z u \
\Psi (x)\\
\end{eqnarray*}
we see that $Z$ is fibrewise the complete intersection of two degree
two equations,
and it does never  intersect the line at infinity $ \{ u=0\}$: whence
we get fibrewise a zero dimensional subscheme $Z_x$ of length $4$.

The fibre subscheme is smooth (iff it consists of $4$ distinct points)
exactly when the Jacobian determinant  $ 4 zw - \Phi \Psi \neq  0.$
This shows that
outside of $R$ we have a local diffeomorphism between $Z$ and $X$.

In the case where $\Phi  \equiv 0 , \Psi   \equiv 0,$ $Z$ is a smooth
(complex) submanifold,
        by the implicit function theorem, if and only if
the two divisors $\{ f=0\}$ and $\{ g=0\}$ are smooth and intersect
transversally.

If $|\Phi| < < 1$, $|\Psi| < < 1$, by uniform continuity on compact sets,
the  real Jacobian matrix has still rank  $4$, thus $Z$ is still a real
submanifold. The equation of the branch divisor $\De$ is obtained
eliminating $z,w$ from the above three equations (see \cite{c-w2}).

\QED

\begin{oss}\label{2}
Assume that $Z$ is a smooth submanifold. Then $R$ is smooth at the
points  $ p \in Z$ where
the fibre $Z_x$ has multiplicity $\mu$ at most two. Also $\De$ is
smooth  at a point
$x$ such that the fibre $Z_x$ consists of three distinct points.
\end{oss}

\Proof
Since $ p \in R$, the multiplicity $ \mu \geq 2$. If equality holds,
there is a linear form $ v (z,w)$ such that $Z$ is locally
defined by $ v^2 =  \varphi(x)$. Hence $Z$ is smooth if and only
if $ R = \{   \varphi(x) = 0\}$ is smooth.

\qed

We want to analyse now the singularities that $\De$ has, for a
general choice of
$\Phi, \Psi$.  In view of the above remark, it suffices to consider points
$x$ such that either

I) there is only one point $ p \in R$ above it, and with multiplicity
$3$ or $4$

II) the fibre $R_x$ consists of two points with multiplicity $2$.

\begin{df}
A point $x$ is said to be a {\bf trivial singularity} of $\De$ if it
is either a point where $ f = \Phi = 0$, or a point where $ g = \Psi = 0$.

$Z$ is said to be {\bf mildly general } if

(*)  at any point where $ f = \Phi = 0$,
we have $ g \cdot\Psi \neq 0$, and symmetrically at any point where
$ g = \Psi = 0$, we have $ f \cdot\Phi \neq 0$

\end{df}

\begin{lem}\label{trivial}
If $Z$ is mildly general, case II) occurs exactly for the trivial singularities
of $\De$. Moreover, case I) occurs only for points  $p \in R$ with
$w(p) z(p) \neq 0$,
and never leads to a point of $Z_x$ of multiplicity $4$.
\end{lem}

\Proof
Observe that over a trivial singularity, for instance with  $ f = \Phi = 0$,
we have $z=0$, and the fibre consists of two distinct points of $R$
if $Z$ is mildly general, since $ g(x) \neq 0$.

Conversely, if we have a point $ p \in R_x$ with $z(p)=0$ and
$w(p) \neq 0$, then $ g(x) \neq 0$, and $\Phi \Psi (x) = 0$
since $ 4 zw -  \Phi \Psi (x) = 0$.

If $\Phi(x)=0$ we have a trivial singularity,
since $ f = z^2 - \Phi w $ also vanishes at $p$.

If instead $\Phi(x) \neq 0$,
then the equation $ f = z^2 - \Phi w $ singles out only one value for $w$ with
$z=0$, and there is no other point $p' \in R_x$. In fact $p' \neq p$,
        $p' \in R_x$ implies $w(p')=0$, contradicting $ g(x) \neq 0$ 
(recall that
$\Psi (x) = 0$!). Thus this case cannot occur.

We argue similarly for the case $w(p)=0$ and
$z(p) \neq 0$.

CLAIM : the case where $w(p)=0= z(p)=0$ cannot occur for a point $ p \in R_x$.

Since otherwise we would have $ f(x)= g(x) = \Phi \Psi (x) = 0$
contradicting the hypothesis
that $Z$ is mildly general.
\qed

We can finally consider the case $w(p) z(p) \neq 0$, hence $\Phi \Psi
(x) \neq 0$.

We can write our two local equations  as

\begin{eqnarray*}
         & w =\frac{1}{\Phi (x)} (z^2 - f(x)) \
\\
         & (z^2 - f(x))^2 -  z  \ \Psi (x) {\Phi (x)}^2 - g(x) {\Phi (x)}^2 = 0.
\\
\end{eqnarray*}
\indent The form of the second equation implies that the sum of the four roots
equals
$0$.

Hence we conclude that not all the four roots are equal, otherwise
we would have $z=0$ as fourtuple root, implying that $f(x) = \Phi (x)
\Psi (x) =  g(x) \Psi (x) =0$,
whence we would have a trivial singularity, contradicting $w(p) z(p) \neq 0$.
The last assertion is then proven.

Assume that there are two roots of multiplicity $2$, thus our
equation would have the
form $ (z-a)^2  (z + a)^2$, i.e. $ (z^2 - a^2 )^2$. Then the coefficient of $z$
would vanish, contradicting $\Phi \Psi (x) \neq 0$.

\QED

\begin{lem}
If $Z$ is mildly general, and  $|\Phi| < < 1$, $|\Psi| < < 1$,
the singular points of $\De$ which are not
the trivial singularities of $\De$, i.e., those coming from case I),
occur only for points  $p \in Z_x$
of multiplicity $3$ which lie over  arbitrarily  small neighbourhoods
of the points with $ f(x) = g(x)
=0$.
\end{lem}

\Proof
By the previous lemma, we are looking for points $x$ where $\Phi \Psi
(x) \neq 0$,
and where the equation
$$(z^2 - f(x))^2 -  z  \ \Psi (x) {\Phi (x)}^2 - g(x) {\Phi (x)}^2 = 0$$
has a triple root $a$.
Then the above monic equation has the form
$$  (z-a)^3  (z + 3 a) = z^4 - 6 a^2 z^2 + 8 a^3 z - 3 a^4.$$

Comparing the coefficients, we obtain
\begin{eqnarray*}
         & - 8 a^3  =  \Psi (x) {\Phi (x)}^2
\\
         & f (x) = 3 a^2.
\\
\end{eqnarray*}
Whence, for $|\Phi| < < 1$, $|\Psi| < < 1$, also $|a| < < 1$,
therefore also $ |f(x)| < < 1$;
by symmetry we also obtain $ |g(x)| < < 1$ and our assertion is proven.

\QED

We want to describe more precisely the singularities of $\De$ in the case where
$X$ is a complex surface.

\begin{prop}
\label{singularities}
Assume that $X$ is a compact complex surface, that the two divisors
$\{ f=0\}$ and $\{ g=0\}$ are
smooth and intersect transversally in a set $M$ of $m$ points.
Then, for $|\Phi| < < 1$, $|\Psi| < < 1$, and for $Z$  mildly general,
$Z$ is a smooth submanifold of
$\VV$, and the projection
$\pi : \VV \ra X$ induces  a finite covering  $ Z \ra X$ of degree
$4$, with smooth orientable ramification divisor $ R : = \{ 4 zw =
\Phi \Psi \}$, and
with branch divisor  $\De $ having as singularities precisely

1) $ 3m$ cusps lying
(in triples) in an arbitrarily small
neighbourhood of $M$, and moreover

2) the trivial singularities, which are nodes if
the curve $\{ f=0\}$
intersect transversally $\{ \Phi = 0\}$, respectively
if the curve  $\{ g=0\}$
intersect transversally $\{ \Psi = 0\}$.
\end{prop}

\Proof
Since we showed that $Z$ is smooth, by remark \ref{2} it follows that
$R$ is smooth, except
possibly for the points of type I).

But in a neighbourhood of a point
with $ f(x) = g (x) = 0$, $R$ is a small differentiable deformation
of the nodal holomorphic  curve  $R_0$ defined ( for $ \Phi (x) 
\equiv 0, \Psi (x)
\equiv 0$) by the
equation
$ 4 zw = 0 , z^2 = f(x) , w^2 = g(x) $.

Working locally at a point where $ f(x) = g (x) = 0$, we may assume that
$ f(x) , g (x) $ are local holomorphic coordinates $(x_1,x_2)$ so that $Z$
is defined by
$$ x_1 = z^2 - w \Phi , x_2 = w^2 -  z  \Psi ,$$
where $z,w$ are local coordinates and $\Phi, \Psi$ do not vanish.

As it was shown in \cite{c-w2}, section 3, in a new system of coordinates
we may assume w.l.o.g. $ \Phi \equiv \Psi \equiv \eta $, where $\eta$ 
is a small
non zero constant.

  Thus $R$ is  smooth, and  it is
orientable being diffeomorphic
to the holomorphic curve  that one gets for $ \Phi (x)  \equiv \eta 
\equiv  \Psi (x)$.

We argue then as
in \cite{c-w2} and we find that $\De$ has exactly three cusps as
singularities for
      $\eta  < < 1$.

For the case (2) of the trivial singularities, it sufffices to observe that
the equation of $\De$ is given by
$$- f^2 g^2 - \frac{9}{8} fg (\Phi \Psi)^2 + (\Psi)^2 f^3 +
        (\Phi)^2 g^3 + \frac{27}{16^2} (\Phi \Psi)^4.$$
At a point where $f=\Phi=0$, and where $g \neq 0, \Psi \neq 0$,
if we assume that $f , \Phi$ give a local diffeomorphism to a
neighbourhood of the origin
in $\C^2$, the equation becomes of order two, and with leading term
$- f^2 g^2(0) +
        (\Phi)^2 g^3 (0)$, thus it defines an ordinary quadratic singularity.

\QED

Since $R$ is orientable, it follows that if the trivial singularities
of $\De$ are nodes,
then the intersection number of the two branches is $\pm 1$; the local
selfintersection number is  exactly equal to $+1$ when the two orientations of
the two branches combine to yield the natural (complex) orientation of $X$.

We are going now to calculate the precise number of singularities, and of the
positive, respectively negative nodes of $\De$ in the special case we
are interested in,
namely, of a perturbed bidouble cover of the quadric.

\begin{df}\label{separately}
Given  four integers $ a,b,c,d\in\N_{\geq3}$,
the so called {\bf dianalytic perturbations  of simple bidouble covers}  are
the {\rm 4}-manifolds
defined by equations
\begin{eqnarray*}
         & z^2= f_{(2a,2b)}(x_0 , x_1; y_0 , y_1) + w \
\Phi _{(2a-c,2b-d)}(x_0 , x_1; y_0 , y_1)\\
         & w^2=g_{(2c,2d)}(x_0 , x_1; y_0 , y_1) + z \
\Psi _{(2c-a,2d-b)}(x_0 , x_1; y_0 , y_1)\\
\end{eqnarray*}
where $f,g$ are bihomogeneous  polynomials of respective bidegrees
$(2a,2b)$,

        $(2c,2d)$,
        and where  $\Phi, \Psi$ are
polynomials in the  ring $$\C[ x_0 , x_1, y_0 , y_1, \bar{x_0 },
\bar{x_1},\bar{y_0 },\bar{y_1}],$$
        bihomogeneous   of
respective bidegrees
$(2a-c,2b-d)$,  $(2c-a,2d-b)$ (the ring is bigraded here by setting the
degree of  $\bar{x_i }$ equal to
$(-1,0)$ and the degree of  $\bar{y_i }$ equal to
$(0, -1)$ ).

We choose moreover $\Phi, \Psi$ to belong to the subspace where
all monomials are either separately holomorphic or antiholomorphic,
i.e., they admit no factor of the form $x_i \bar{x_j }$ or
of the form $y_i \bar{y_j }$.
\end{df}

\begin{oss}\label{antihol}
Indeed we shall pick $\Phi$ (respectively : $\Psi$) to be a product
$$\Phi = \Phi_1(x)  \Phi_2 (y),$$ where $ \Phi_1$ is a product
of
linear forms, either all holomorphic or all antiholomorphic,
and similarly for  $\Phi_2$ (resp. : for $\Psi_1$, $\Psi_2$).

Under this condition, at a point where $\Phi$ vanishes simply,
we get $$\Phi = \mbox{ unit }\cdot \Phi_1,$$
or $\Phi = \mbox{ unit }\cdot \Phi_2$.

\end{oss}

As explained in \cite{catcime}(pages 134-135) to an antiholomorphic
homogeneous polynomial
$ P(\bar{x_0 },
\bar{x_1})$ of degree $m$ we associate a differentiable section $p$ of
the tensor power $\LL^{\otimes m}$ of the tautological (negative) subbundle,
such that $ p (x) : = P( 1, \bar{x})$ inside a big disc $B(0,r)$
in the complex line
having centre at the origin and  radius $r$.

We shall make the following assumptions on the polynomials $f,g, \Phi, \Psi$:

\begin{enumerate}\label{hyp}
\item
The algebraic curves $C : = \{ f = 0 \}$ and $D : = \{ g = 0 \}$ are smooth
\item
$C$ and $D$ intersect transversally at a finite set $M$ contained in
the product $B(0,r)^2$  and at
these points both curves have non vertical tangents
\item
for both curves $C$ and $D$ the first projection on $\PP^1$ is a
simple covering
and moreover the vertical tangents of $C$ and $D$ are all distinct
\item
the associated perturbation is sufficiently small and mildly general
\item
the trivial singularities are nodes with non vertical tangents,
contained in $B(0,r)^2$
\item
the cusps of $\De$ have non vertical tangent
\end{enumerate}

Consider  the branch curve $\De$ of the perturbed bidouble cover.
Then the defining equation $\de$ is bihomogeneous of bidegree $(4
(a+c), 4 (b+d))$.

We proved that $\De$ has exactly $ k : = 12 (ad + bc)$ cusps, coming from the
        $ m : = 4 (ad + bc)$ points of the set $ M = C \cap D$.

        If $\Phi, \Psi$
are holomorphic, then $\De$ has exactly
$$\nu : = 4 ( 2 ab + 2 cd - ad - bc)$$
(positive) nodes.

In general, this number $\nu$ equals
the difference $\nu^+ - \nu^-$ between the number of positive and
the number of negative nodes.

In fact the nodes of $\De$ occur only for the trivial singularities.

A trivial singularity $ f = \Phi = 0$ yields a node with  the same tangent cone
  as $$ f^2 - \Phi^2 \cdot g(0).$$

By remark \ref{antihol} we may assume that there are local holomorphic
coordinates $f, \hat{\Phi} $ such that either $\Phi = \mbox{ unit 
}\cdot \hat{\Phi}$
or $\Phi = \mbox{ unit }\cdot \overline {\hat{\Phi}}$.

In the first case we get the tangent cone of a holomorphic node, hence
a {\bf positive} node, in the second case we
get the tangent cone of an antiholomorphic node, hence
a {\bf negative} node.

\begin{lem}
The number $t_f$ of vertical tangents for the curve $C$ is
$  t_f = 4 ( 2 ab - a) $, while the  number $t_g$ of vertical
tangents for the curve $D$ is
$  t_g = 4 ( 2 cd - c) $. For a general choice of $\Phi, \Psi$, the
number of vertical tangents
$t$ of $\De$ equals $ t = 2 t_f + 2 t_g + m$.
The genus of the Riemann surface $R$ equals $ g =  1 + 16 (a+c)(b+d)
- 4 (a+b+c+d) - k -\nu $,
where $ k : =  3 m = 12 (ad + bc)$ is the number of cusps of $\De$ and
$\nu : = \nu^+ - \nu^-$ is the number of nodes of $\De$, counted with sign.
\end{lem}

\Proof
The curve $C$ has bidegree $(2a,2b)$ thus its canonical divisor has degree
$ 2 g(C) - 2 = 2a (2b-2) + 2b (2a-2)= 8 ab - 4 a - 4b$. The first
projection has degree $2b$,
thus by the Hurwitz' formula the number of ramification points is $ 2
g(C) - 2 + 4 b$.
Since the morphism is simple, the number $t_f$ of vertical tangent
equals the number of ramification points. Similarly for $D$.

In the case where $\Phi \equiv 0,  \Psi \equiv 0,$ the ramification
curve $R_0$ has $m$ double
points, and is composed of two double covers of $C \cup D$, branched
on the set of double points.

Therefore, if $g$ is the  genus of $R$, $ 2g - 2 = 2  (2 g(C) - 2) +
2  (2 g(D) - 2) + 4m =
16 ab - 8 a - 8 b + 16 cd  - 8c - 8 d + 16 (ad + bc) = 16 (a+c)(b+d)
- 8 (a+b+c+d) $.

Now, the first projection yields a map  $p : R \ra\PP^1_{\C}$ of
degree $4 (b+d)$,
thus the formula for $t$ is again an application of Hurwitz' formula
since $p|_R$ is finite (the map from $\De$ to $\PP^1_{\C}$ being finite, for
$\Phi, \Psi$ sufficiently small).

$R$ is smooth orientable, and $p$ has only a finite number of critical
values. Thus $p$ is oriented, and indeed a ramified covering. For
general $\Phi, \Psi,$ it
has simple ramification, with distinct critical values. Hence the
number of critical
values
is given by $ 2g -2 + 8 (b+d) = 16 (a+c)(b+d) - 8 (a+c) = 2 t_f + 2
t_g + 16 (ad + bc)=
2 t_f + 2 t_g +  4m = 2 t_f + 2 t_g + m + k$.

We conclude observing that the number of critical values equals the
number $k$ of
(non vertical) cusps plus the number $t$ of vertical tangents.

\QED

\begin{oss}
1) Observe that the fact that $p : R \ra\PP^1_{\C}$ has a finite
number of critical points implies
immediately that $p$ is finite, since there is then a finite set in
$\PP^1_{\C}$
such that over its complement we have a covering space. Thus $p$ is
orientation preserving
and each non critical point contributes positively to the degree of
the map $p$.

2) Consider now a real analytic 1-parameter family $Z(\eta)$ such that
the perturbation terms
$\Phi_{\eta}, \Psi_{\eta} \ra 0$. Then we have a family of ramification points
$P_i({\eta})$ which tend
to the ramification points of $p_0 : R_0 \ra\PP^1_{\C}$. For $R_0$, a
double cover of
$C \cup D$, we have $2t_f$ ramification points lying over the $t_f$
vertical tangents
of $C$,  $2t_g$ ramification points lying over the $t_g$ vertical tangents
of $D$, whereas the $m$ nodes of $R_0$ are limits of $4$ ramification points,
three cusps and a simple vertical tangent. This is the geometric reason why
$t =  2 t_f + 2 t_g + m$; it also tells us where to look for the
vertical tangents of $\De$.
\end{oss}

In the next section we shall also consider the curves $C$ and $D$ as
small real deformations
of reducible real curves with only nodes (with real tangents) as singularities:
according to the classical notation, these will be called {\bf improper nodes}.

Moreover, we shall also consider real polynomials $\Phi, \Psi,$ such that the
trivial singularities are also given by real points. By the results of
\cite{c-w2}, from the perturbation of each proper node of $C \cup D$,
i.e., of a point of $ M = C \cap D$, we shall obtain  a real vertical  tangent,
a real cusp, and two immaginary cusps.

%%%%%%%%%%%%%%%%%%%%%%%%%%%%%%%%%
\section{the basic model}

\newcommand{\nodalC}{C^\times}
\newcommand{\nodalD}{D^\times}
\newcommand{\nodalf}{f^\times}
\newcommand{\nodalg}{g^\times}

In this section we will provide a model with some choices of
$f,g,\Phi,\Psi$ to which the reasoning of the previous chapter
applies. In particular they are explicit enough to collect
a lot of geometric information and to determine the braid
monodromy factorization in the sequel.

Let us begin with real polynomials $\nodalf$ and $\nodalg$
defining nodal curves $\nodalC$, respectively $\nodalD$.

For the sake of simplicity we replace here bihomogeneous
polynomials $ f(x_0, x_1 ; y_0, y_1 )$ by their restrictions to the
affine open set
$x_0 =1 , y_0 =1$ and write $ f(x,y)$ for $ f(1, x ; 1, y )$.

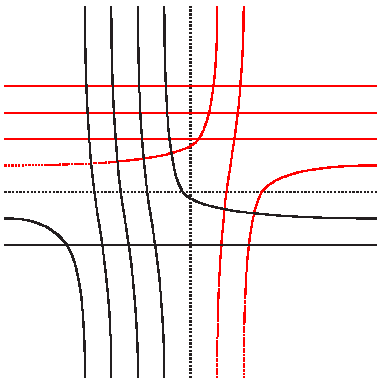
\begin{figure}
\centering
\begin{picture}(150,150)(-75,-75)
\bezier{100}(-70,0)(0,0)(70,0)
\bezier{100}(0,-70)(0,0)(0,70)

%f%
\rot{
\bezier{200}(-70,20)(0,20)(70,20)
\bezier{200}(-70,30)(0,30)(70,30)
\bezier{200}(-70,40)(0,40)(70,40)

\bezier{100}(15,10)(20,40)(20,70)
\bezier{100}(15,10)(10,-20)(10,-70)

\bezier{100}(3,20)(10,30)(10,70)
\bezier{100}(3,20)(-4,10)(-70,10)
\bezier{100}(27,0)(20,-10)(20,-70)
\bezier{100}(27,0)(34,10)(70,10)
}
%g%

\bezier{400}(-70,-20)(0,-20)(70,-20)

\bezier{200}(-15,-10)(-10,-40)(-10,-70)
\bezier{200}(-15,-10)(-20,20)(-20,70)
\bezier{200}(-25,-10)(-20,-40)(-20,-70)
\bezier{200}(-25,-10)(-30,20)(-30,70)
\bezier{200}(-35,-10)(-30,-40)(-30,-70)
\bezier{200}(-35,-10)(-40,20)(-40,70)

\bezier{200}(-47,-20)(-40,-30)(-40,-70)
\bezier{200}(-47,-20)(-54,-10)(-70,-10)
\bezier{200}(-3,0)(-10,10)(-10,70)
\bezier{200}(-3,0)(4,-10)(70,-10)

\end{picture}
     \caption{Real part of $\Delta_f$ and $\Delta_g$ for
       $(a,b,c,d)=(1,2,2,1)$}
\label{delta f and g}
\end{figure}

We let $\nodalf : =f_1\cdots f_{2b}$ defining  $f_i : = y-2i$
     for i=$2, \dots , 2b$ but setting
$$
f_1: =(y-2)\prod_{i=1}^{2a}(x-2i)+\prod_{i=1}^{2a-1}(x-2i-1),
$$
and we define similarly $\nodalg : =g_1\cdots g_{2d}$ setting $g_i=y+2i$  for
i=$2, \dots , 2d$ except that we set
$$
g_1: =(y+2)\prod_{i=1}^{2c}(x+2i)-\prod_{i=1}^{2c-1}(x+2i+1).
$$

\bigskip
\begin{oss}
\label{halfplane}
Note that the equations of $f_1$ and $g_1$ are chosen in such a way
that their zero sets are graphs of rational functions $\tilde{f}_1$,
respectively
     $\tilde{g}_1$, of
$x$ which, regarded as maps from
$\C$ to
$\C$, preserve the real line. Moreover, $\tilde{f}_1$ preserves both the upper
and lower halfplane, while $\tilde{g}_1$ exchanges them.

In fact, $\tilde{f}_1$ has no critical points on $\PP^1_{\R}$ and
$$ deg (\tilde{f}_1 |_{ \PP^1_{\R} }) \colon \PP^1_{\R }\ra \PP^1_{\R}  =
2a = deg (\tilde{f}_1 ).$$
Similarly for $\tilde{g}_1$.
\end{oss}

Next we define a real analytic 1-parameter family of dianalytic perturbation
terms
$\Phi_\eta$ and $\Psi_\eta$ introducing the following notation:
$$
\begin{array}{cclclcl}
\Phi_{(2a-c,0)} & = &
\multicolumn{4}{l}{
\left\{ \begin{array}{ll}
\displaystyle
\prod_{j=1}^{2a-c}(x-4a-2j), & \text{ if } 2a\geq c,\\[2mm]
\displaystyle
\prod_{j=1}^{c-2a}(\bar x-4a-2j), & \text{ if } 2a<c.
\end{array}\right.}\\
\Phi_{(0,2b-d)} & = &
\multicolumn{4}{l}{
\left\{ \begin{array}{ll}
\displaystyle
\prod_{j=1}^{2b-d}(y-4b-2j), & \text{ if } 2b\geq d,\\[2mm]
\displaystyle
\prod_{j=1}^{d-2b}(\bar y-4b-2j), & \text{ if } 2b<d.
\end{array}\right.}\\
\Psi_{(2c-a,0)} & = &
\multicolumn{4}{l}{
\left\{ \begin{array}{ll}
\displaystyle
\prod_{j=1}^{2c-a}(x+4c+2j), & \text{ if } 2c\geq a,\\[2mm]
\displaystyle
\prod_{j=1}^{a-2c}(\bar x+4c+2j), & \text{ if } 2c<a.
\end{array}\right.}\\
\Psi_{(0,2d-b)} & = &
\multicolumn{4}{l}{
\left\{ \begin{array}{ll}
\displaystyle
\prod_{j=1}^{2d-b}(y+4d+2j), & \text{ if } 2d\geq b,\\[2mm]
\displaystyle
\prod_{j=1}^{b-2d}(\bar x+4d+2j), & \text{ if } 2d<b.
\end{array}\right.}\\[15mm]
\Phi_\eta & : =&  \eta\Phi_{(2a-c,2b-d)} & = &
\eta\Phi_{(2a-c,0)}\Phi_{(0,2b-d)}\\[4mm]
\Psi_\eta & : =&  \eta\Psi_{(2c-a,2d-b)} & = &
\eta\Psi_{(2c-a,0)}\Psi_{(0,2d-b)}
\end{array}
$$

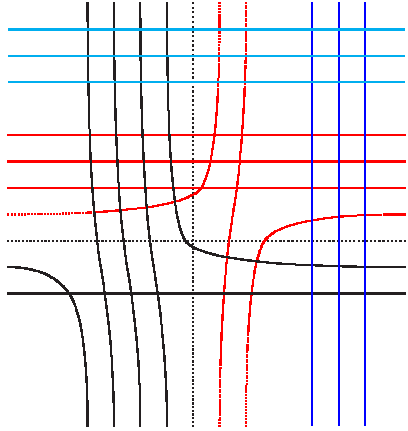
\begin{figure}
\centering
\begin{picture}(170,170)(-75,-75)
\bezier{100}(-70,0)(0,0)(80,0)
\bezier{100}(0,-70)(0,0)(0,90)

%f%
\rot{
\bezier{200}(-70,20)(0,20)(80,20)
\bezier{200}(-70,30)(0,30)(80,30)
\bezier{200}(-70,40)(0,40)(80,40)

\bezier{100}(15,10)(20,40)(20,90)
\bezier{100}(15,10)(10,-20)(10,-70)

\bezier{100}(3,20)(10,30)(10,90)
\bezier{100}(3,20)(-4,10)(-70,10)
\bezier{100}(27,0)(20,-10)(20,-70)
\bezier{100}(27,0)(34,10)(80,10)
}
%g%

\bezier{400}(-70,-20)(0,-20)(80,-20)

\bezier{200}(-15,-10)(-10,-40)(-10,-70)
\bezier{200}(-15,-10)(-20,20)(-20,90)
\bezier{200}(-25,-10)(-20,-40)(-20,-70)
\bezier{200}(-25,-10)(-30,20)(-30,90)
\bezier{200}(-35,-10)(-30,-40)(-30,-70)
\bezier{200}(-35,-10)(-40,20)(-40,90)

\bezier{200}(-47,-20)(-40,-30)(-40,-70)
\bezier{200}(-47,-20)(-54,-10)(-70,-10)
\bezier{200}(-3,0)(-10,10)(-10,90)
\bezier{200}(-3,0)(4,-10)(80,-10)

%\put(22,-16){$0$}
{\color{cyan}
{\color{blue}
\put(45,-70){\line(0,1){160}}
\put(55,-70){\line(0,1){160}}
\put(65,-70){\line(0,1){160}}
}
{
\put(-70,60){\line(1,0){150}}
\put(-70,70){\line(1,0){150}}
\put(-70,80){\line(1,0){150}}
}
}
\end{picture} \caption{Zero sets of $\Phi$ (cyan) and $\Psi$
(magenta) added}
\label{phipsi}
\end{figure}

Given $\nodalf$, consider the polynomial $f : \nodalf + c_f$,
where $c_f$ is a small constant; likewise    consider  $g : \nodalg + c_g$.
Adding these small  constants to the
respective equations of
$\nodalC$ and
$\nodalD$
we get polynomials $f$ and $g$ which define
smooth curves $C$ and $D$. We have more precisely:

\begin{prop}
\label{hypprop}
If the constants $c_f$, $c_g$, $\eta$
are chosen sufficiently small, the polynomials $f, g, \Phi, \Psi$ thus obtained
satisfy  the following list of hypotheses,
cf. \pageref{hyp},  for $r >8(a+b+c+d)$:
\begin{enumerate}\label{hyp2}
\item
The algebraic curves $C = \{ f = 0 \}$ and $D = \{ g = 0 \}$
are smooth
\item
$C$ and $D$ intersect transversally at a finite set $M$ contained in
the product $B(0,r)^2$  and at
these points both curves have non vertical tangents
\item
for both curves $C$ and $D$ the first projection on $\PP^1$ is a
simple covering
and moreover the vertical tangents of $C$ and $D$ are all distinct
\item
the associated perturbation is mildly general
\item
the trivial singularities of $\Delta$
are nodes with non vertical tangents,
contained in $B(0,r)^2$
\item
the cusps of $\De$ have non vertical tangent
\end{enumerate}
\end{prop}

\proof
Property (1) follows obviously from Bertini's theorem, since we have pencils
without base points.

Property (2) follows since the conditions are open
conditions which hold for the nodal curve $\nodalC$ and $\nodalD$, as we
shall show in  lemma \ref{hypnodal} below.

By the same lemma the critical values for $\nodalC\cup\nodalD$ are
distinct, hence property (3) is proven if we show that the two vertical
tangents arising from each improper node
do not map to the same critical value.
If this were the case, the local degree of the projection would be 4 and not 2.

To prove property (4) that the associated perturbation is mildly general,
it suffices  -- by definition and continuity -- to note
that no three of the functions
$\nodalf,\nodalg,\Phi,\Psi$ have a common zero.

With prop.\ \ref{singularities} we deduce from (4) that the trivial
singularities are nodes, which yields the first part of (5).
Moreover the local analysis in the proof
of prop.\ \ref{singularities} shows that the two tangents at the node
are close to the tangent to $C\cup D$ at that point.
So the second part of (5) follows by continuity, since the points
of $\nodalC\cup \nodalD$ with vertical tangents
are not on the zero sets of $\Phi$ and $\Psi$.

(6) As shown in \cite{c-w2},section 3,  and already used in proposition
\ref{singularities},  we may assume to have
that  $f, g$ are local coordinates $U,V$ , and that $\Phi, \Psi$ are
small constants $\eta_1 , \eta_2$,
so that the local equation of $\De$ is
$$ -U^2 V^2 - \frac{9}{8} UV c^2 + c^2 V^3 + U^3 + \frac{27}{16^2} c^4 ,$$
$c$ being the ratio $ \frac{\eta_1}{\eta_2} = \frac{\Phi}{\Psi} $.

Following the proof in loc. cit., we set $ c = \lambda^3$, and see that
setting $ U = \lambda^4 u_0$,$ V = \lambda^4 v_0$ the equation becomes
$$ -u_0^2 v_0^2 - \frac{9}{8} u_0 v_0 c^2 + c^2 v_0^3 + u_0^3 +
\frac{27}{16^2} .$$

Then the directions
  of the tangent at
the three cusps in the coordinates $(U,V)$ are non constant functions of $c$.

Hence for a general choice of the ratio $c$ between $\phi(0)$ and
$\Psi(0)$
these tangents are not vertical.

\qed

\begin{lemma}
\label{hypnodal}
The critical values of the vertical projection of
$\nodalC\cup\nodalD$ are all distinct and real. Moreover
\begin{enumerate}
\item[(*)]
$\nodalC$ and $\nodalD$ intersect transversally at a finite set
$M^\times$ contained in
the bidisk $B(0,r)^2 $  and at
these points both curves have non vertical tangents
\end{enumerate}
\end{lemma}

\proof
Each critical point is a singular point of the union
$\nodalC\cup\nodalD$ since the components of these curves are
horizontal sections.
Given two components of $\nodalC\cup\nodalD$ we can easily check that the
number of  real  points in their intersection equals the intersection number
of these two components: hence  all singularities are real, they are nodes,
and moreover all
     critical values are real.

Recall now that each component of $C^\times$ except $\{ f_1 = 0 \}$
is a horizontal line, and similarly for $\nodalD$ we have only the
component $\{ g_1 = 0 \}$ which is not a horizontal line.
Hence the nodes are of three types: type $C$, where $\{ f_1 = 0 \}$
meets a horizontal component of $\nodalD$, mixed type,
where $ f_1 = g_1 = 0 $, and type $D$.

A vertical line through a point of mixed type cannot contain points of
type $C$, nor $D$, so we have to exclude that a vertical line contains
a point of type $C$ and a point of type $D$. But we simply observe that
a point of type $C$ has positive abscissa, while a point of type $D$
has negative abscissa.

To prove (*) we note again that all singularities are nodes and
we check from the equations that their coordinates are
bounded by $8(a+b+c+d)$. Since all components of
$\nodalC\cup\nodalD$ map biregularly to the base,
no component has a vertical tangent.
\qed

\vspace*{1cm}

We introduce now a geometric system of paths
in the $\C$-line  with coordinate $x$ (this is
the  part not at infinity of the target space
$\PP^1$ of the vertical projection):
$$
\cg_{\nu,\rho},\quad (\nu,\rho)\in
\mathcal{V}=\mathcal{V}^+\cup\mathcal{V}^-
$$
where
$$
\mathcal{V}^+ =\bigg\{ (\nu,\rho)\,\Big|\;
1\leq \phantom{-}\nu\leq 4a,
\Big\{\begin{array}{rll}
     1<\rho\leq \hspace{1.7mm}2b & \text{if }\nu \text{ odd}\\
-2d\leq\rho\leq -1 & \text{if }  \nu \text{ even}
\end{array} \bigg\}
,$$
respectively
$$
\mathcal{V}^- =\bigg\{ (\nu,\rho)\,\Big|\;
1\leq -\nu\leq 4c,
\Big\{\begin{array}{rll}
     1\leq\rho\leq \hspace{1.7mm}2b & \text{if }\nu \text{ odd}\\
-2d\leq\rho< -1 & \text{if }  \nu \text{ even}
\end{array} \bigg\}
$$
correspond to the singular points of $\nodalC\cup\nodalD$:
$\mathcal{V}^+$ corresponds to those with positive abscissa,
$\mathcal{V}^-$ corresponds instead to those with negative abscissa.

Observe in fact that the number of positive critical values
is $2a(2b+2d-1)$ and thus equals  the cardinality of $\mathcal{V}^+$.

$\mathcal{V}^+$, respectively $\mathcal{V}^-$, is ordered
lexicographically according to $|\nu|$ and
$\rho$, and these orderings correspond each to the
total order of the critical values on the real line
given by the absolute value.

The total order of $\mathcal{V}$ is then determined by requiring
     that elements in $\mathcal{V^+}$ preceed those in $\mathcal{V^-}$.

The $\cg_{\nu,\rho}$'s  are a system of simple paths based
at the origin, following each other in counterclockwise order,
each $\cg_{\nu,\rho}$  going around
(again in a  counterclockwise way) the corresponding critical value.

In fact we can choose the paths with
$\nu$ positive to lie in the upper half plane
except in a neighbourhood of the critical value
(lying on the positive real line) which they encircle.
Similarly the paths with $\nu$ negative can be
chosen to lie in the lower half plane except
     in a neighbourhood of the critical value
(lying on the negative real line) which they encircle.

\begin{figure}
\begin{picture}(340,70)(-20,-5)

\bezier{600}(-10,30)(150,30)(310,30)
\bezier{200}(150,0)(150,30)(150,60)

\multiput(10,30)(5,0){4}{\circle*{2}}
\multiput(45,30)(5,0){4}{\circle*{2}}
\multiput(80,30)(5,0){4}{\circle*{2}}
\multiput(115,30)(5,0){4}{\circle*{2}}
\multiput(0,30)(35,0){4}{\circle*{2}}

\multiput(170,30)(12,0){3}{\circle*{2}}
\multiput(214,30)(12,0){2}{\circle*{2}}
\multiput(250,30)(12,0){3}{\circle*{2}}
\multiput(294,30)(12,0){2}{\circle*{2}}

\put(170,30){\circle{8}}
\put(306,30){\circle{8}}

\bezier{100}(150,30)(160,36)(166,30)
\bezier{1000}(150,30)(210,90)(302,30)

\multiput(165,34)(0,3){3}{\circle*{1}}

\put(130,30){\circle{6}}
\put(0,30){\circle{8}}

\bezier{100}(150,30)(140,25)(133,30)
\bezier{1000}(150,30)(90,-30)(4,30)

\multiput(135,26)(0,-3){3}{\circle*{1}}

\end{picture} \caption{The geometric system for the curve of figure
\ref{delta f and g}}
\label{gamma}
\end{figure}
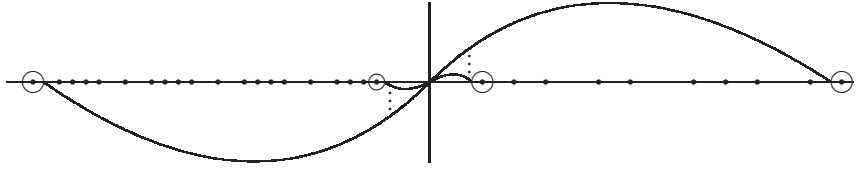

To explain the choice of the index set in more detail let us observe
some restrictions for the $y$-coordinate function on
     the different components of
$\nodalC\cup\nodalD$ when the $x$-coordinate is restricted to be
real and positive:
it is constant of value $2i$, resp.\ $-2j$ for the components
defined by $f_i$ resp.\ $g_j$ with $i,j>1$
and bounded between $-1$ and $-2$ for the component
defined by $g_1$.
It is real or infinity for the component $C_1$ defined by $f_1$ and
is monotone increasing outside of the poles.

Since the $y$-coordinate on $C_1$ takes value in $[2,3]$ for
$x = 0$,
we can precisely describe in which order the other components
are intersected by $C_1$ as  $x$ increases.

First $C_1$ intersects the other components $C_2, \dots C_{2b}$
of $\nodalC$
in this order at points labelled by $(1,2), \dots (1,2b)$;
     then the $x$-coordinate
crosses a pole of $f_1$ and  $y$ changes sign.

Next $C_1$ intersects the components $D_{2d}, \dots D_{1}$
of $\nodalD$ in this order at points  labelled by $(2,-2d), \dots (2,-1)$.

After
that the $y$ coordinate changes sign again, and the same pattern of
intersections  repeats itself over and over till the $y$-coordinate
finally tends to
$2$ as the
$x$-coordinate tends to infinity.

We can finally summarize  how our choice of the index set relates to
this sequence of intersections:
\begin{itemize}
\item
a critical value with $\nu>0$ corresponds to an intersection
with $C_1$,
\item
each time the $y$-coordinate changes sign, $\nu$ is
incremented by one,
\item
a critical value with $\rho$ positive corresponds to an intersection
with the component defined by $f_\rho$,
\item
a critical value with $\rho$ negative corresponds to an intersection
with the component defined by $g_{-\rho}$.
\end{itemize}

Of course the situation over the negative real line can be described
similarly and we have:
\begin{itemize}
\item
a critical value with $\nu$ negative corresponds to an intersection
point with the component defined by $g_1$,
\item
each time the $y$-coordinate changes sign over a negative
real $x$, the absolute value of the index $\nu$ is incremented
by one.
\end{itemize}

In the reference fibre over the origin the roots of $\nodalf$ are
real and positive,
in fact the root of $f_i$ is $2i$ for $i>1$ and between $2$ and $3$
for $i=1$.
Similarly the roots of $\nodalg$ are real and negative, the root of $g_j$
is $-2j$ except for $g_1$ when the root is between $-2$ and $-1$.

We define now arcs in the reference fibre, contained in the upper half-plane
and connecting two roots according to the following labelling
$$
\tau_{i,j}:\quad
\left\{
\begin{array}{ll}
i=1, j>1& \text{between roots  of } f_1,f_j\\
i>0, j<0& \text{between roots  of } f_i,g_{|j|}\\
i<0,j=-1& \text{between roots  of } g_1,g_{|i|}
\end{array}
\right.
$$

and subject to the following rule: for each fixed $i$ (resp.  $j$)
the arcs $\tau_{i,j}$
follow each other in counterclockwise order without intersecting
outside of the end points,
and otherwise two such arcs  outside of the end points  intersect
transversally at at most one point.

These arcs are thus determined up to an isotopy fixing a neighbourhood of
infinity.

\begin{prop}
The braid monodromy factorization for the curve
$\nodalC \cup \nodalD $ with respect to the geometric basis
$\cg_{\nu,\rho}$ can be given in terms of the mapping classes
$\b_{i,j}$, the half-twists associated to the arcs $\tau_{i,j}$,
in the following form:
$$
\cg_{\nu,\rho} \quad\mapsto\quad
\left\{
\begin{array}{ll}
\b_{1,\rho}^2 & \text{if } \nu>0,\\
\b_{\rho,-1}^2 & \text{if } \nu<0
\end{array}\right.
$$
\end{prop}

\proof
Since all critical points on $\nodalC\cup\nodalD$
for the vertical projection are nodes, the mapping class
associated to any element $\cg_{\nu,\rho}$ of the
geometric system is necessarily a full twist on some arc
between two roots.

These two roots are in fact the endpoints of $\tau_{1,\rho}$,
resp.\ $\tau_{\rho,-1}$, -- depending on the sign of $\nu$ --
due to the way the index pair $(\nu,\rho)$ remembers
the two components which intersect at the critical value
encircled by $\cg_{\nu,\rho}$.

It thus remains to prove that the mapping class associated
to $\cg_{\nu,\rho}$ can be obtained from an arc in the upper
half plane.

Let us consider the case $\nu>0$. Then we specify the geometric loop
$\cg_{\nu,\rho}$ as follows: a simple
    arc in the upper half plane, which will be called the 'tail',
connecting the base point with
a point on the real line, which we shall call the 'fork point',
followed by a circle around the critical value corresponding to
$(\nu,\rho)$ run
once counterclockwise before the loop goes back
along the tail.

When $\cg_{\nu,\rho}$ runs the circle
around the critical value we get a mapping class in the fibre over the
fork point of $\cg_{\nu,\rho}$
which is obviously the full twist on the straight real arc
connecting the roots associated to the two components of $\nodalC \cup \nodalD$
which intersect at the critical point.

Along the tail of $\cg_{\nu,\rho}$ all roots stay fixed for the horizontal
components, the root associated to $g_1=0$ moves
in the lower half plane, the root associated to $f_1=0$ moves
in the upper half plane, cf.\ remark \ref{halfplane}.

We may thus associate continuously
to each point $x$ on the tail of $\cg_{\nu,\rho}$ an arc $\tau_x$
in the fibre $F_x$ over $x$,
such that the mapping class associated to $\cg_{\nu,\rho}$ truncated at $x$
is given by the full twist on $\tau_x$.

In fact this arc can be chosen to be the straight
arc except for $x$ on $\cg_{\nu,\rho}$ close to the origin, because
for $x=0$ we would get the straight arc through other roots.

By continuity also the arc associated to $x=0$ may be chosen to lie
in the upper half plane and it is thus isotopic to $\tau_{1,\rho}$,
as we wanted to show.

The argument in the case $\nu<0$ is the same.

\qed

We have now reached a point where we
can  draw the first two consequences concerning
the braid monodromy homomorphism associated to the
vertical projection of the branch curve $\Delta$.

By continuity we may deduce first the following claim about the
mapping classes associated to the divisor $\Delta$
along the paths $\cg_{\nu,\rho}$.
To state it properly
let $A_{i,j}$ be respective closed neighbourhoods of the arcs $\tau_{i,j}$,
each containing, among the roots  of $\nodalf\nodalg$,
only the two roots
which form the endpoints of $\tau_{i,j}$.

\begin{prop}
\label{cont}
Let $\Delta$ be defined as in proposition \ref{hypprop} and choose
sufficiently small
constants so that it moreover lies in a sufficiently close neighbourhood of
$\nodalC\cup\nodalD$. Then the following holds true:
\begin{enumerate}
\item
each $\cg_{\nu,\rho}$ encircles all critical values
of $ p_1 \colon\De \ra \PP^1$ which tend
to the corresponding critical value of $\nodalC\cup\nodalD$.
\item
the mapping class associated to $\cg_{\nu,\rho}$
is supported on $A_{i,j}$ if $\b_{i,j}^2$ is the mapping class
associated in the case of the nodal curve $\nodalC\cup\nodalD$.
\end{enumerate}
\end{prop}

Of course we will henceforth tacitly assume to be in the situation
of the above proposition.

The other consequence is that all mapping classes associated to the paths
$\cg_{\nu,\rho}$ are given by braids lying in a distinguished
subgroup of the braid group.
To define this subgroup consider the $2$-cable homomorphism
$\br_n$ to $\br_{2n}$, which associates to each braid
a corresponding braid with doubled strands, eg.
\begin{center}
\begin{picture}(300,50)
\bezier{100}(45,10)(45,20)(40,25)
\bezier{100}(40,25)(35,30)(35,40)
\bezier{100}(35,10)(35,18)(38,22)
\bezier{100}(42,28)(45,32)(45,40)
\bezier{80}(54,10)(54,25)(54,40)
\bezier{80}(63,10)(63,25)(63,40)

\put(83,23){$\mapsto$}

\bezier{100}(130,10)(130,20)(120,25)
\bezier{100}(120,25)(110,30)(110,40)
\bezier{100}(132,10)(132,21)(122,26)
\bezier{100}(122,26)(112,31)(112,40)

\bezier{100}(110,10)(110,19)(117,23)
\bezier{100}(124,28)(130,32)(130,40)
\bezier{100}(112,10)(112,18)(118,22)
\bezier{100}(125,27)(132,31)(132,40)

\bezier{80}(140,10)(140,25)(140,40)
\bezier{80}(142,10)(142,25)(142,40)
\bezier{80}(150,10)(150,25)(150,40)
\bezier{80}(152,10)(152,25)(152,40)
\end{picture}
\end{center}

\newcommand{\cbr}{\ON{CBr}}
\newcommand{\fbr}{\ON{CBr}}

\subsection{Definition}
Define the subgroup $\cbr_n\subset\br_{2n}$ of 2-
cable braids to be the
subgroup generated by the image of  the $2$-cable homomorphism,
and by the half-twist inside each cable (i.e.,\ each doubled strand).
\\

The half-twists inside the cables generate a free abelian
normal subgroup $\cong \Z^n$
leading, via   the $2$-cable homomorphism, to a representation of $\cbr_n$ as
the semidirect product
$$
\cbr_n \cong \Z^n >\!\!\!\lhd \br_n .
$$

Here $ \Z^n : = \langle \sigma_{1 + 2i} \rangle_{i=0, \dots,  n-1}$ and
$ \br_n : = \langle \sigma'_{j} \rangle_{j=1, \dots , n-1}$,
where $$  \sigma'_{j} : =  \sigma_{ 2j}\sigma_{2j -1 }\sigma_{2j + 
1}\sigma_{2j}.$$

\begin{lemma}
\label{framed}
The braid monodromy associated to $\Delta$
maps the homotopy class
of any path $\cg_{\nu,\rho}$
to the subgroup $\fbr_n$ of $\br_{2n}$.

Moreover  the braid monodromy
associated to $\nodalC\cup\nodalD$ along $\cg_{\nu,\rho}$
is the image of the above through the natural
projection $\fbr_n\to\br_n$.
\end{lemma}

\proof
Braids in $\fbr$ can be realized as braided ribbons, where each
doubled strand is considered as the  boundary of a twisted
ribbon. Along any path $\cg_{\nu,\rho}$ the braid associated to
$\Delta$ can be extendend to a ribbon: we simply connect
each pair of roots regenerating from a single point of $C\cup D$
by a straight arc and observe that such arcs never intersect,
since the roots of $C$ and $D$ stay sufficiently away from each other
along $\cg_{\nu,\rho}$.

The second assertion is obvious.

\qed

In the next step we investigate the local braid monodromy factorizations
in more detail for the regenerations of proper and
improper nodes of $\nodalC\cup\nodalD$.

In the  case of a proper node we get a cusp-cluster, four critical
points of which
three are cusps and the last is a vertical tangency point.
In the case of an improper node we get a cluster of vertical tangents, four
critical points, all of which are vertical tangency points.

The braid monodromy for the regeneration into a cusp-cluster
has been thoroughly investigated before in \cite{c-w2},
where the braid monodromy factorization type  has been
determined up to
Hurwitz equivalence and simultaneous conjugation.
\footnote{  The factorization was explicitly described for a particular
choice of coordinates $(x,y)$.

However, any other coordinate change, $\hat{x} = F(x)$,
$\hat{y} = G(x, y)$, induces another description of the fibre
bundle  induced by $ p (x,y) := x$ on the complement of $\De$.
This results into a Hurwitz equivalence (action of the base diffeomorphism $F$)
followed by a simultaneous conjugation  (action  of the diffeomorphism
$G$ on the reference fibre).

}

We take this investigation over again in order to determine the
    Hurwitz equivalence class of the factorization
from the datum of  the product of the
factors.

\newcommand{\braid}{\cg}

\begin{prop}
\label{hur1}
The braid monodromy factorization of a (regenerated) cusp-cluster with
product $\s_2^3\s_1\s_3\s_2\s_1^2\s_3^2$ is
Hurwitz-equivalent to the factorization
$$
\s_2^3 \circ  \s_1 \s_3 \s_2 \s_3\inv \s_1\inv \circ \s_1^3 \circ \s_3^3.
$$
\end{prop}

\proof
By the results of  \cite{c-w2} the two factorizations
are equivalent up to Hurwitz-equivalence and
simultaneous conjugation.
Since these two operations commute
the braid monodromy factorization of our cusp-cluster
is Hurwitz equivalent to a factorization
$$
(\s_2^3)^\braid
\circ (\s_1 \s_3 \s_2 \s_3\inv \s_1\inv)^\braid
\circ (\s_1^3)^\braid
\circ (\s_3^3)^\braid
$$
where, as usual,  $\be^\braid : = \braid \be \braid^{-1}.$

Since Hurwitz equivalence does not change the product of a
factorization, it follows that
    $\braid$ is in the centralizer of the braid
$\s_2^3\s_1\s_3\s_2\s_1^2\s_3^2$.

    From the work of Gonzales-Meneses and Wiest \cite{gmw} we may
read off two elements,
$\Delta_4 : =\s_1\s_2\s_3\s_1\s_2\s_1$ and $\s_1\inv\s_2\s_1$,
which generate this centralizer.

To prove our claim it thus suffices to give a Hurwitz equivalence
between the factorization given in the above formula and
its conjugates by $\s_1\inv\s_2\s_1$, respectively $\Delta_4$.
Below we first give expressions for the two respective conjugated
factorizations using the fact that $\s_i^{\Delta_4}=\s_{4-i}$.
\begin{eqnarray*}
&&
(\s_2^3)^{\s_1\inv\s_2\s_1}
\circ (\s_1 \s_3 \s_2 \s_3\inv \s_1\inv)^{\s_1\inv\s_2\s_1}
\circ (\s_1^3)^{\s_1\inv\s_2\s_1}
\circ (\s_3^3)^{\s_1\inv\s_2\s_1}\\[2mm]
& = &
(\s_2^3)^{\s_2\s_1\s_2\inv}
\circ (\s_1 \s_2\inv \s_3 \s_2 \s_1\inv)^{\s_1\inv\s_2\s_1}
\circ (\s_1^3)^{\s_1\inv\s_2}
\circ (\s_3^3)^{\s_1\inv\s_2}\\[2mm]
& = &
(\s_2^3)^{\s_2\s_1}
\circ ( \s_3 )^{\s_1\inv\s_2\s_1\s_1 \s_2\inv}
\circ \s_1\inv\s_2\s_1^3\s_2\inv\s_1
\circ (\s_3^3)^{\s_1\inv\s_2}\\[2mm]
& = &
\s_1^3
\circ ( \s_3 )^{\s_1\inv\s_1\inv\s_2\s_2 \s_1}
\circ \s_1\invv\s_2^3\s_1^2
\circ (\s_3^3)^{\s_1\inv\s_2}\\[2mm]
& = &
\s_1^3
\circ ( \s_3 )^{\s_1\invv\s_2^2}
\circ (\s_2^3)^{ \s_1\invv}
\circ (\s_2^3)^{\s_1\inv\s_3\inv}
\\[5mm]
&  &
(\s_2^3)^{\Delta_4}
\circ (\s_1 \s_3 \s_2 \s_3\inv \s_1\inv)^{\Delta_4}
\circ (\s_1^3)^{\Delta_4}
\circ (\s_3^3)^{\Delta_4}\\[2mm]
& = &
\s_2^3 \circ \s_3 \s_1 \s_2 \s_1\inv \s_3\inv
\circ \s_3^3 \circ \s_1^3\\
\end{eqnarray*}
Finally we show that these two elements are in fact in the same
Hurwitz-equivalence class by performing elementary Hurwitz
operations and rewriting of the factors.
Here the notation $\sim_{(i\,i+1)}$ is employed to
denote a simple Hurwitz move affecting the positions $i$
and $i+1$.

\begin{align*}
     & \quad \s_2^3
\circ \s_1 \s_3 \s_2 \s_3\inv \s_1\inv
\circ \s_1^3
\circ \s_3^3\\
     \sim_{(23)} & \quad
\s_2^3
\circ \s_1^3
\circ (\s_1 \s_3 \s_2 \s_3\inv \s_1\inv)^{\s_1^{-3}}
\circ \s_3^3\\
     = & \quad
\s_2^3
\circ \s_1^3
\circ \s_1\invv \s_3 \s_2 \s_3\inv \s_1^2
\circ \s_3^3\\
     \sim_{(12)} & \quad
\s_1^3
\circ (\s_2^3)^{\s_1^{-3}}
\circ \s_1\invv \s_3 \s_2 \s_3\inv \s_1^2
\circ \s_3^3\\
     \sim_{(23)} & \quad
\s_1^3
\circ \s_1\invv \s_3 \s_2 \s_3\inv \s_1^2
\circ (\s_2^3)^{\s_1\invv \s_3 \s_2\inv \s_3\inv \s_1^2\s_1^{-3}}
\circ \s_3^3\\
     = & \quad
\s_1^3
\circ \s_1\invv \s_3 \s_2 \s_3\inv \s_1^2
\circ (\s_2^3)^{\s_3 \s_1\invv \s_2\inv \s_3\inv \s_1\inv}
\circ \s_3^3\\
     \sim_{(34)} & \quad
\s_1^3
\circ \s_1\invv \s_3 \s_2 \s_3\inv \s_1^2
\circ \s_3^3
\circ (\s_2^3)^{\s_3^{-3}\s_3 \s_1\invv \s_2\inv \s_3\inv \s_1\inv}\\
     = & \quad
\s_1^3
\circ \s_1\invv \s_3 \s_2 \s_3\inv \s_1^2
\circ \s_3^3
\circ (\s_1^3)^{\s_3\invv \s_1\invv \s_2\inv \s_3\inv \s_2}\\
     = & \quad
\s_1^3
\circ \s_1\invv \s_3 \s_2 \s_3\inv \s_1^2
\circ \s_3^3
\circ (\s_1^3)^{\s_3\invv \s_1\invv \s_3 \s_2\inv \s_3\inv}\\
     = & \quad
\s_1^3
\circ \s_1\invv \s_3 \s_2 \s_3\inv \s_1^2
\circ \s_3^3
\circ (\s_1^3)^{\s_3\inv \s_1\inv \s_1\inv \s_2\inv}\\
     = & \quad
\s_1^3
\circ \s_1\invv \s_3 \s_2 \s_3\inv \s_1^2
\circ \s_3^3
\circ (\s_2^3)^{\s_1\inv \s_3\inv}\\
     = & \quad
\s_1^3
\circ \s_1\invv \s_3 \s_2 \s_3\inv \s_1^2
\circ \s_3^3
\circ (\s_3^3)^{\s_1\inv \s_2}\\
     \sim_{(23)} & \quad
\s_1^3
\circ (\s_3^3)^{\s_1\invv \s_3 \s_2 \s_3\inv \s_1^2}
\circ \s_1\invv \s_3 \s_2 \s_3\inv \s_1^2
\circ (\s_2^3)^{\s_3\inv \s_1\inv}\\
     = & \quad
\s_1^3
\circ (\s_3^3)^{\s_1\invv \s_3 \s_2}
\circ \s_1\invv \s_3 \s_2 \s_3\inv \s_1^2
\circ (\s_2^3)^{\s_3\inv \s_1\inv}\\
     = & \quad
\s_1^3
\circ (\s_2^3)^{\s_1\invv}
\circ \s_1\invv \s_3 \s_2 \s_3\inv \s_1^2
\circ (\s_2^3)^{\s_3\inv \s_1\inv}\\
     = & \quad
\s_1^3
\circ \s_1\invv \s_2^3\s_1^2
\circ (\s_2)^{ \s_1\invv \s_3}
\circ (\s_2^3)^{\s_3\inv \s_1\inv}\\
     \sim_{(23)} & \quad
\s_1^3
\circ (\s_2)^{\s_1\invv \s_2^3\s_1^2 \s_1\invv \s_3}
\circ \s_1\invv \s_2^3\s_1^2
\circ (\s_2^3)^{\s_3\inv \s_1\inv}\\
     = & \quad
\s_1^3
\circ (\s_2)^{\s_1\invv \s_2^2\s_2 \s_3}
\circ (\s_2^3)^{\s_1\invv}
\circ (\s_2^3)^{\s_3\inv \s_1\inv}\\
     = & \quad
\s_1^3
\circ (\s_3)^{\s_1\invv \s_2^2}
\circ (\s_2^3)^{\s_1\invv}
\circ (\s_2^3)^{\s_1\inv \s_3\inv}
\\[5mm]
& \quad
\s_2^3
\circ \s_1 \s_3 \s_2 \s_3\inv \s_1\inv
\circ \s_1^3
\circ \s_3^3\\
     \sim_{(34)}& \quad
\s_2^3
\circ \s_1 \s_3 \s_2 \s_3\inv \s_1\inv
\circ \s_3^3
\circ \s_1^3\\
\end{align*}
Comparison now shows that we got the same elements.
\qed

In the case of a regeneration into a cluster of vertical tangents
we have to understand first the braid monodromy factorization type
up to Hurwitz equivalence and simultaneous conjugation.

\begin{prop}
\label{vt-fact}
The braid monodromy factorization associated to a cluster of
tangents corresponding to an improper node  is given by
\[
\s_2\s_3\s_2\inv \circ \s_1\s_2\s_1\inv \circ
\s_2\s_3\s_2\inv \circ \s_1\s_2\s_1\inv
\]
up to Hurwitz equivalence and simultaneous conjugation.
\end{prop}

\proof
Without loss of generality, we assume that
we have an improper node of $\nodalC$.

We may take a coordinate change not altering
the $x$ coordinate such that the local equation of
    $\nodalC$  is $u(x,y) (y^2-x^2 )= 0$,
where $u(x,y)  \neq 0$ (we also have $\nodalg(x,y)  \neq 0$).

Consider first a small perturbation with $\Psi=0$:
then the branch curve is given by
\[
{\nodalf}^2{\nodalg}^2-{\nodalg}^3\Phi^2={\nodalg}^2({\nodalf}^2-{\nodalg}\Phi^2)=0.
\]
The vertical tangents can be located from the non-unit factor.
\[
\frac{\partial}{\partial y}({\nodalf}^2-{\nodalg}\Phi^2)=0
\Leftrightarrow \frac{\partial}{\partial y}((y^2-x^2)^2-{\nodalg}
u^{-2} \Phi^2)=0
\Leftrightarrow
\]
\[
\quad \Leftrightarrow\quad 2y(y^2-x^2)=\frac{\partial}{\partial y}
({\nodalg}u^{-2}\Phi^2)
\]
To simplify the equations we take a local model
based on an apt  choice of $\Phi$
    $({\nodalg}u^{-2}\Phi^2) \equiv \epsilon^2$,
where $\epsilon$ is a small constant,  which makes the right hand side vanish.
The two equations have solutions
$y=0$ and $x^4= \epsilon^2$,
which leads to four solutions close to $(0,0)$.

Then it is not hard to see that the same braid monodromy
factorization is obtained from the model
\[
(y^2-x^2-\epsilon)(y^2-x^2+\epsilon)=0
\]
for which the factorization is as claimed.

The actual perturbation $\Delta$ of $\nodalC\cup\nodalD$
is linked to the above through a path of perturbations
which constantly have four vertical tangency points at
the cluster.
So the local fibration does not change and
therefore also the equivalence class of braid monodromy
factorization for Hurwitz equivalence and
simultaneous conjugation stays the same.

\qed

\begin{prop}
\label{hur2}
Suppose that a braid monodromy factorization of a cluster of tangents
is  Hurwitz plus  conjugation  equivalent to the factorization given by
$$
(\s_1\inv\s_2\s_1 ) \circ  (\s_2\inv\s_3\s_2 ) \circ
(\s_1\inv\s_2\s_1) \circ ( \s_2\inv\s_3\s_2).
$$
Then the two factorizations are equivalent via  Hurwitz equivalence
and simultaneous conjugation by a power of $\s_1$
if and only if  their products are equal.
\end{prop}

\proof
Again we deduce immediately
that  their products are equal iff the first
factorization is Hurwitz equivalent to the given one
via conjugation by some braid centralizing the product.

The product  is however $(\Delta_4\s_1^{-2}\s_3^{-2})^2$ and
by \cite{gmw} its centralizer is generated by $\Delta_4$ and $\s_1$.

This immediately shows that the given equivalence preserves the product.

In the other direction, in order to prove our claim it thus suffices 
to give a Hurwitz equivalence
between the factorization of the claim and
its conjugate by $\Delta_4$.

We first give the expression for this conjugate
factorization using that $\Delta_4\s_i\Delta_4\inv=\s_{4-i}$.
\begin{eqnarray*}
& &
(\s_1\inv\s_2\s_1)^{\Delta_4} \circ
(\s_2\inv\s_3\s_2)^{\Delta_4}  \circ
(\s_1\inv\s_2\s_1)^{\Delta_4} \circ
(\s_2\inv\s_3\s_2)^{\Delta_4}  \\[2mm]
& = &
(\s_3\inv\s_2\s_3) \circ
( \s_2\inv\s_1\s_2) \circ
(\s_3\inv\s_2\s_3) \circ
(\s_2\inv\s_1\s_2)\\[2mm]
\end{eqnarray*}
We finish providing an elementary chain of transformations, obtained 
using only Hurwitz equivalence, simultaneous
conjugation by $\s_1^2$ and rewriting of factors, which leads from 
the factorization of the claim
to its conjugate by $\Delta_4$.

\begin{eqnarray*}
& &
(\s_1\inv\s_2\s_1)^{\s_1^2} \circ
(\s_2\inv\s_3\s_2)^{\s_1^2} \circ
(\s_1\inv\s_2\s_1)^{\s_1^2} \circ
(\s_2\inv\s_3\s_2)^{\s_1^2} \\[2mm]
& = &
(\s_1\s_2\s_1\inv) \circ
(\s_2\inv\s_3\s_2)^{\s_1^2} \circ
(\s_1\inv\s_2\s_1)^{\s_1^2} \circ
(\s_2\inv\s_3\s_2)^{\s_1^2} \\[2mm]
&  \sim_{(12)(34)} &
(\s_2\inv\s_3\s_2)^{\s_1\s_2\s_1\inv\s_1^2} \circ
(\s_1\s_2\s_1\inv) \circ
(\s_2\inv\s_3\s_2)^{\s_1\s_2\s_1\inv\s_1^2} \circ
(\s_1\s_2\s_1\inv) \\[2mm]
&  = &
(\s_2\inv\s_3\s_2)^{\s_1\s_2\s_1} \circ
(\s_1\s_2\s_1\inv) \circ
(\s_2\inv\s_3\s_2)^{\s_1\s_2\s_1} \circ
(\s_1\s_2\s_1\inv) \\[2mm]
&  = &
(\s_2\inv\s_3\s_2)^{\s_2\s_1\s_2} \circ
(\s_1\s_2\s_1\inv) \circ
(\s_2\inv\s_3\s_2)^{\s_2\s_1\s_2} \circ
(\s_1\s_2\s_1\inv)\\[2mm]
&  = &
(\s_3)^{\s_2\s_1} \circ
(\s_1\s_2\s_1\inv) \circ
(\s_3)^{\s_2\s_1} \circ
(\s_1\s_2\s_1\inv)  \\[2mm]
&  = &
(\s_3)^{\s_2} \circ
(\s_1\s_2\s_1\inv) \circ
(\s_3)^{\s_2} \circ
(\s_1\s_2\s_1\inv) \\[2mm]
&  = &
(\s_2\s_3\s_2\inv) \circ
(\s_1\s_2\s_1\inv) \circ
(\s_2\s_3\s_2\inv) \circ
(\s_1\s_2\s_1\inv) \\[2mm]
&  = &
(\s_3\inv\s_2\s_3) \circ
(\s_2\inv\s_1\s_2) \circ
(\s_3\inv\s_2\s_3) \circ
(\s_2\inv\s_1\s_2)
\end{eqnarray*}
\qed

\begin{cor}
\label{admiss}
In the situation of the above proposition the two factorizations are equivalent
up to Hurwitz equivalence and creation/deletion of admissible pairs.
\end{cor}

\proof
It suffices to show that simultaneous conjugation by $\s_1^2$
can be induced by Hurwitz moves and creation/deletion of admissible
pairs. This can be done as follows:
\begin{align*}
& \quad \s_1\inv\s_2\s_1 \circ \s_2\inv\s_3\s_2 \circ
\s_1\inv\s_2\s_1 \circ \s_2\inv\s_3\s_2\\
\sim_c & \quad  \s_1^{-2} \circ \s_1^{2} \circ
\s_1\inv\s_2\s_1 \circ \s_2\inv\s_3\s_2 \circ
\s_1\inv\s_2\s_1 \circ \s_2\inv\s_3\s_2
\\
\sim & \quad \s_1^{-2} \circ
(\s_1\inv\s_2\s_1)^{\s_1^2} \circ
(\s_2\inv\s_3\s_2)^{\s_1^2} \circ
(\s_1\inv\s_2\s_1)^{\s_1^2} \circ
(\s_2\inv\s_3\s_2)^{\s_1^2}
\circ \s_1^2\\
\sim & \quad \s_1^{-2} \circ \s_1^2 \circ
(\s_1\inv\s_2\s_1)^{\s_1^2} \circ
(\s_2\inv\s_3\s_2)^{\s_1^2} \circ
(\s_1\inv\s_2\s_1)^{\s_1^2} \circ
(\s_2\inv\s_3\s_2)^{\s_1^2}\\
\sim_d & \quad
(\s_1\inv\s_2\s_1)^{\s_1^2} \circ
(\s_2\inv\s_3\s_2)^{\s_1^2} \circ
(\s_1\inv\s_2\s_1)^{\s_1^2} \circ
(\s_2\inv\s_3\s_2)^{\s_1^2}.
\end{align*}

A braid monodromy factorization is associated to the branch
curve $\Delta$ once we fix  a geometric system of paths
for the critical values of $\Delta$ in the base
   and a system of arcs between the roots in the
reference fibre $F_0$.

Concerning the geometric system of paths, we shall construct them
later as a `regeneration'
of
the old system of paths $\cg_{\nu,\rho}$: ie.\ given a path $\cg_{\nu,\rho}$
we shall determine a sequence of consecutive paths -- actually four of them --
in the new system, having the same tail, and such that their product
is homotopic to
$\cg_{\nu,\rho}$. In particular, the choice for these paths will be a local
one around the critical values for $\nodalC \cup \nodalD$.

The choice of the local factorization  of each $\cg_{\nu,\rho}$ corresponds
to a choice inside the Hurwitz equivalence class of the braid factorization
whose product is  the braid associated to $\cg_{\nu,\rho}$.

We shall now concentrate on the systems of arcs.

%

%%%%%%%%%%%%%%%%%%%%%%%%%%%%%

%%%%%%%%%%%%%%%%%%%%%%%%%%%%%%%%%
\section{proof of the factorization theorem}

\begin{prop}
\label{c-mono}
Up to homotopy in small discs containing pairs $D'_i,D''_i$ or $B_\jmath',
B_\jmath''$ the punctured fibre $F_0$ is given by the following picture
\\
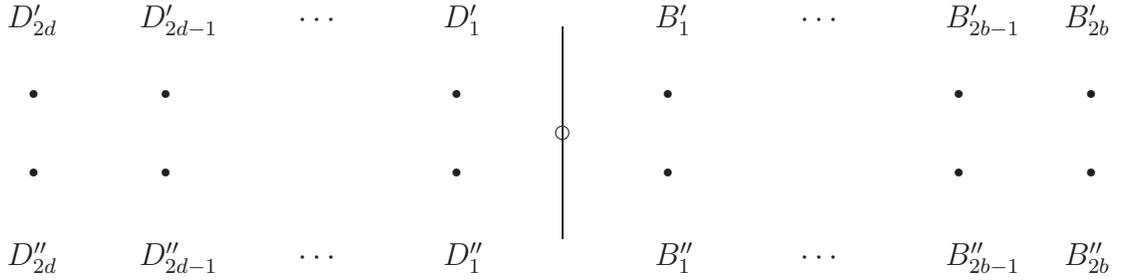
\begin{picture}(420,120)(-210,-60)

\put(0,0){\circle{5}}
\put(0,-40){\line(0,1){80}}

\put(40,15){\circle*{3}}
\put(150,15){\circle*{3}}
\put(200,15){\circle*{3}}
\put(-40,15){\circle*{3}}
\put(-150,15){\circle*{3}}
\put(-200,15){\circle*{3}}
\put(40,-15){\circle*{3}}
\put(150,-15){\circle*{3}}
\put(200,-15){\circle*{3}}
\put(-40,-15){\circle*{3}}
\put(-150,-15){\circle*{3}}
\put(-200,-15){\circle*{3}}

\put(-210,-50){$D''_{2d}$}
\put(-160,-50){$D''_{2d-1}$}
\put(-45,-50){$D''_{1}$}
\put(35,-50){$B''_{1}$}
\put(145,-50){$B''_{2b-1}$}
\put(190,-50){$B''_{2b}$}
\put(-210,40){$D'_{2d}$}
\put(-160,40){$D'_{2d-1}$}
\put(-45,40){$D'_{1}$}
\put(35,40){$B'_{1}$}
\put(145,40){$B'_{2b-1}$}
\put(190,40){$B'_{2b}$}

\put(-100,-50){$\cdots$}
\put(-100,40){$\cdots$}
\put(90,-50){$\cdots$}
\put(90,40){$\cdots$}

\end{picture}\\
and the covering monodromy $\theta$ of the perturbed bidouble cover,
with respect to the origin,
is given by
$$
\begin{array}{rcl@{\hspace{2cm}}rcl}
\oo_{D'_j} & \mapsto & (12), & \oo_{B'_j} & \mapsto & (13),\\
\oo_{D''_j} & \mapsto & (34), &
\oo_{B''_j} & \mapsto & (24),\\
\end{array}
$$
where $\oo_P$ is any simple closed path around $P$
not crossing the imaginary axis $i\RR\subset F_0$.
\end{prop}

\proof
As we remarked before, in the unperturbed Galois cover case we have
real points $D_i$,
$B_\jmath$, only, the  $D_i$'s with positive real coordinate,
the  $B_j$'s with negative real coordinate.s
And the cover monodromy is $(12)(34)$, resp.\ $(13)(24)$
for paths which do not cross the imaginary axis.

The deformation then splits each branch point into two, hence the
corresponding monodromies (which are transpositions)
must be $(12)$ and $(34)$, resp.\ $(13)$
and
$(24)$ (which commute).
After a suitable homotopy the points are in the positions  given by
the picture.

\qed

Referring to  the  above figure we introduce the following notation for
a system of arcs, which
are uniquely determined (up to homotopy) by their endpoints and the
property that
they are monotonous in the real coordinate and do not pass below
any puncture; i.e.,  if they share the real coordinate with a puncture
they have larger imaginary coordinate.
$$
\begin{array}{ll@{\hspace*{2cm}}ll}
p_{i}: & B'_i, B''_i, &
q_{\imath}: & D'_\imath, D''_\imath,\\
a_{ij}: & B'_i, B'_j, &
b_{\imath\jmath}: & D'_ \imath, D'_\jmath, \\
c_{ij}: & B''_i, B''_j, &
d_{\imath\jmath}: & D''_\imath, D''_\jmath, \\[2mm]
u'_{i\jmath}: & B'_i, D'_\jmath &
u''_{i\jmath}: & B''_i, D''_\jmath
\end{array}
$$

\begin{oss}
We can safely assume that any such arc is fully contained in one of
the previously defined closed
sets $A_{i,j}$ (neighbourhood of the arc $\tau_{i,j}$).
\end{oss}

\begin{prop}
\label{v-cluster}
The braid monodromy factorization associated to $\Delta$ and to
a subsystem of paths
refining $\cg_{\nu,\rho}$ is given by
\begin{enumerate}
\item
\qquad
$\s_{a_{1i}}\circ \s_{p_1}^2\s_{c_{1i}}\s_{p_1}^{-2}\circ
\s_{a_{1i}}\circ \s_{p_1}^2\s_{c_{1i}}\s_{p_1}^{-2}$,
\qquad
if $\nu>0$ and $i=\rho>0$,
\item
\qquad
$\s_{b_{1j}}\circ \s_{q_1}^2\s_{d_{1j}}\s_{q_1}^{-2}\circ
\s_{b_{1j}}\circ \s_{q_1}^2\s_{d_{1j}}\s_{q_1}^{-2}$,
\qquad
if $\nu<0$ and $j=-\rho>0$,
\end{enumerate}

up to Hurwitz equivalence and simultaneous conjugation
by some power of $\s_{p_1}^2$ in case (1), resp.\ by  $\s_{q_1}^2$
in case (2).
\end{prop}

\proof
The proofs of both cases are the same modulo an appropriate
exchange of indices, so we consider the first case only.

Since the paths refine $\cg_{\nu,\rho}$ the associated mapping
classes are supported on $A_{1,i}$ for $i=\rho$, by $(3)$ of prop.\
\ref{cont}.

Hence by prop.\ \ref{vt-fact} the factorization coincides
with the factorization of the claim up to simultaneous
conjugation by a mapping class (of this disc) supported on $A_{1,i}$.
The products (i.e., the product of the factorization 1. , and the 
product of the
factorization we are looking for, yielding the braid monodromy associated to
$\cg_{\nu,\rho}$)
are therefore  equal up to that conjugation.

To proceed we have to get more information about the second product:
being associated to the path $\cg_{\nu,\rho}$ it  belongs
to the subgroup $\ON{CBr}_n$ of $\br_{2n}$ by lemma \ref{framed}.

Next, the product belongs to the subgroup of braids
which are supported on $A_{1,i}$.
This subgroup is generated by $\s_{a_{1i}},\s_{c_{1i}},\s_{p_{1}}, 
\s_{p_{i}}$ and
we can easily determine a presentation of its intersection
with $\ON{CBr}_n$:
$$
\langle \s_{p_{1}}, \s_{p_{i}}, \Delta_4:=\s_{a_{1i}} \s_{p_{1}}^2 
\s_{c_{1i}}  \s_{p_i}^2\,|\,
\s_{p_1} \s_{p_i} = \s_{p_i} \s_{p_1}, \Delta_4 \s_{p_1} = \s_{p_i} 
\Delta_4, \Delta_4 \s_{p_i} = \s_{p_1}
\Delta_4\rangle
$$
Again by lemma \ref{framed} we know the image of the product under
the quotient map $\ON{CBr}_n\to\br_{2n}$ to be a full twist.
Since both $\s_{p_1}$ and $\s_{p_i}$ belong to the kernel of that map we
deduce now that this product can be written as $\Delta_4^2 \s_{p_1}^k 
\s_{p_i}^{k'}$
with $k+k'=-8$ to match the total degree, which is $4$.

The final fact to be exploited is that the two products are conjugate
inside the group generated by $\s_{a_{1i}},\s_{c_{1i}},\s_{p_1}, \s_{p_i}$.

The first product can be shown to be equal to $\Delta_4^2 
\s_{p_1}^{-4} {\s_{p_i}}^{-4}$
since the two elements $\s_{p_1}$ and $\s_{p_i}$ commute, while
conjugation by $\Delta_4$ exchanges them.
Since moreover the element
$\Delta_4^2$ is central we get conjugation equivalences
$$
\Delta_4^2 \s_{p_1}^k \s_{p_i}^{k'} \sim \Delta_4^2 \s_{p_1}^{-4} 
{\s_{p_i}}^{-4}
\quad\implies\quad
\s_{p_1}^k \s_{p_i}^{k'} \sim  \s_{p_1}^{-4} {\s_{p_i}}^{-4}
$$

But closing the braid on both sides we get a four component link for 
the right hand side,
which  consists of two unlinked copies of the $(4,2)$ torus link.

This soon implies that, since the braids are conjugate,  $k$ and $k'$ 
are both even, hence the
left hand side consists of the  unlinked union of a $(k,2)$ torus link
with a $(k',2)$ torus link.

The non zero linking numbers on the right are equal to $2$, while
the non zero linking numbers on the left are equal to $ 
\frac{|k|}{2}, \frac{|k'|}{2} $.

  Again since the braids are conjugate and the linking numbers
are invariant under conjugation we conclude from $|k|= |k'| = 2$ and $k+k'= -4$
that $k=k'=-2$, so both
products are equal.

So finally we may
invoke prop.\ \ref{hur2} and we are left only
with the ambiguity of a conjugation by a power of $\s_{p_1}$.

This power must however have an even exponent, since the group of 
liftable braids is
left invariant by conjugation by $\s_{p_1}^2$, while  conjugation by $\s_{p_1}$
carries  $\s_{a_{1i}}$ outside of the subgroup of liftable braids.

\qed

\begin{prop}
\label{c-cluster}
Suppose $\nu$ and $\rho$ have opposite sign.
Then the braid monodromy factorization associated to $\Delta$ and to
a subsystem of paths
refining $\cg_{\nu,\rho}$ is given by
$$
\s_{u_{i\jmath}}^3 \circ
\s_{s_{i\jmath}} \circ
\s_{u_{i\jmath}'}^3 \circ
\s_{u_{i\jmath}''}^3
\qquad\text{ where}\quad
\left\{\begin{array}{ll}
i=1,j=-\rho & \text{if } \nu>0,\\
i=\rho, j=1 & \text{if } \nu<0.
\end{array}\right.
$$
up to Hurwitz equivalence.
In particular there is a distinguished refining subsystem of paths
such that the associated factorization is precisely the one given above.
\end{prop}

\proof
Since the paths refine $\cg_{\nu,\rho}$, by prop.\ \ref{cont} (3)
the factors are supported on $A_{i,j}$,
which is a topological disc containing four roots,
$B_i',B_i'',D_j',D_j''$.
Locally the divisor $\Delta$ is a regeneration of
a node between the two components $C$ and $D$, thus
by the results of Catanese and Wajnryb \cite{c-w2} the factorization
coincides with the factorization of the claim up to Hurwitz equivalence and
simultaneous conjugation.
Hence also their products coincide up to conjugation.

These products are both inside $\fbr_2$ and map to the same
braid $\b_{i,j}^2$ under the map $\fbr_2\to\br_4$.
Hence they coincide up to conjugation by elements in the
subgroup generated by the half-twists on the arcs $p_i$ and $q_j$.
Since these half-twists are central in $\fbr_2$,
the products do in fact coincide.

Therefore we may invoke
prop. \ref{hur1} to deduce that the factorization is
Hurwitz equivalent to the one given in the claim.

In particular the Hurwitz action on the factors is induced by
the Hurwitz action on the subsystem of paths, thus we may
pick the subsystem of paths in the Hurwitz orbit which
maps to the given factorization.

\qed

\begin{thm}\label{BMF}
The factorization is a product of factorizations
$$
\begin{matrix}
\underbrace{
\bigg(\b_f \circ
\underbrace{(\s_{p_1}^{\pm2} \circ \cdots \circ
\s_{p_1}^{\pm2})}_{|2b-d|\text{ factors}}
\circ \b_{fg} \bigg)\circ \cdots}_{2a\text{ repetitions}}
\circ
\underbrace{(\s_{p_1}^{\pm2} \circ \cdots \circ
\s_{p_{2b}}^{\pm2})\circ \cdots}_{|2a-c|\text{ repetitions}}
\\[2mm]
\hspace*{30mm}\circ
\underbrace{(\s_{q_1}^{\pm2} \circ \cdots \circ
\s_{q_{2d}}^{\pm2})\circ \cdots}_{|2c-a|\text{ repetitions}}
\circ
\underbrace{
\bigg(\b_g \circ
\underbrace{(\s_{q_1}^{\pm2} \circ \cdots \circ
\s_{q_1}^{\pm2})}_{|2d-b|\text{ factors}}
\circ \b_{gf} \bigg)\circ \cdots}_{2c\text{ repetitions}}
\end{matrix}
$$
where the sign of the exponents $2$ is constant inside a pair of
brackets and is the sign of the number which determines the number
of factors, ie.\ $(2b-d)$ resp.\ $(2a-c)$, $(2c-a)$ or $(2d-b)$,
and
where the $\b$'s further decompose as products of factorizations
\[\begin{matrix}
\b_f & = & \b_{f,2}\circ \cdots \circ \b_{f,2b},\quad
\b_{fg} & = & \b_{fg,2d}\circ \cdots \circ \b_{fg,1},\\
\b_{g} & = & \b_{g,2} \circ \cdots \circ \b_{g,2d},\quad
\b_{gf} & = & \b_{gf,2b} \circ \cdots \circ \b_{gf,1},
      \end{matrix}\]
based on elementary factorizations each having four factors
\[\begin{matrix}
\b_{f,i} & = & \s_{a_{1i}}\circ \s_{p_1}^2\s_{c_{1i}}\s_{p_1}^{-2}\circ
\s_{a_{1i}}\circ \s_{p_1}^2\s_{c_{1i}}\s_{p_1}^{-2},\qquad
\b_{fg,j} & = & \s_{u_{1j}}^3\circ \s_{s_{1j}}\circ
\s_{u_{1j}'}^3\circ \s_{u_{1j}''}^3,\\
\b_{g,j} & = & \s_{b_{1j}}\circ \s_{q_1}^2\s_{d_{1j}}\s_{q_1}^{-2}\circ
\s_{b_{1j}}\circ \s_{q_1}^2\s_{d_{1j}}\s_{q_1}^{-2},\qquad
\b_{gf,i} & = & \s_{u_{i1}}^3\circ \s_{s_{i1}}\circ
\s_{u_{i1}'}^3\circ \s_{u_{i1}''}^3.
      \end{matrix}\]
The elementary factorizations originate in the regeneration of nodes
of the branch curve of the corresponding bidouble Galois-cover.
\end{thm}

\proof
Consider first the set of critical values for the vertical projection
of $\Delta$. They are either arbitrarily close to the critical values
of the vertical projection of $\nodalC\cup\nodalD$ or they are images
of the nodes of $\Delta$ and thus are real.

We choose an associated system of paths in the  base such that:
\begin{enumerate}
\item
the paths associated to positive critical values of nodes belong to the
upper half plane,
\item
the paths associated to negative critical values of nodes belong to the
lower half plane,
\item
the paths associated to the four critical values which regenerate
from a node of $\nodalC\cup\nodalD$ refine the path
$\cg_{\nu,\rho}$ associated to the value of that node.
\end{enumerate}
We have still the choices of elements in the Hurwitz orbits of
subsystems refining each path $\cg_{\nu,\rho}$.
   We use these choices to adjust the system of paths using
prop. \ref{v-cluster} and prop. \ref{c-cluster} and we succeed to get
the elementary
factorizations, at least up to the conjugations mentioned in loc. cit.

The mapping classes associated to the other paths
are obviously  full twists.
We have to handle the problem that there may be critical values
corresponding to different
critical points: this happens if and only if there are nodes coming
from vertical components of $\Phi,\Psi$.

We can handle this difficulty  adding a constant to $\Phi$, resp. $\Psi$.
This resolves the
problem except for the case that the divisor is not ample, ie.,
if $2b=d$ or $2b=d$.

In this case we take a local differentiable perturbation at the nodes.

\qed

We pass now to consider the conjugacy classes of the braids appearing in the
braid monodromy factorization, in order to prove Theorem \ref{bmf}.

\begin{lemma}
\label{admissible}
All full-twists of type $p$ and $q$ are admissible.
\end{lemma}

\proof
The arc of type $p_i$ is homotopic to the union of two paths
$\cg_{D_i'},\cg_{D_i''}$ in the
positive half plane. By prop.\ \ref{c-mono}
we see that the monodromies around the two ends are disjoint
transpositions, hence the
full-twist on $p_i$ is admissible.
The same argument yields the claim for $q_\jmath$.

\qed

%%%
\subsection{Generators of mapping class groups}

\begin{lemma}
\label{triangle-rel}
If $a_{ij},a_{jk},a_{ik}$ are the sides of a topological triangle taken in
positive order,
then
$$
\s_{a_{ij}}\s_{a_{jk}}\quad=\quad
\s_{a_{jk}}\s_{a_{ik}}
$$
\end{lemma}

\proof
This is easily seen from
the identity $\s_1\s_2=\s_2(\s_2\inv\s_1\s_2)$
in the braid group
$\br_3$ which maps homomorphically
onto the group generated
by
$\s_{a_{ij}},\s_{a_{jk}},\s_{a_{ik}}$ via
$$
\s_1 \: \mapsto \:
\s_{a_{ij}},\qquad
\s_2 \: \mapsto \:
\s_{a_{jk}},\qquad
\s_2\inv\s_1\s_2 \: \mapsto \:
\s_{a_{ik}}.
$$
\qed

\begin{lemma}
\label{braid-gen}
The following
subgroups of the braid group $\br_{4(b+d)}$
are identical:
$$
\langle
\s_{a_{1,2}}, \s_{a_{2,3}},\cdots,\s_{a_{n-1,n}} \rangle,
\langle
\s_{a_{1,2}},\cdots,\s_{a_{1,n}} \rangle,
\langle \s_{a_{ij}}, 1\leq
i,j \leq n  \rangle,
$$
\end{lemma}

\proof
By the relation in lemma
\ref{triangle-rel} a group containing
$\s_{a_{ij}},\s_{a_{jk}}$ also
contains $\s_{a_{ik}}$ for any
triple $i,j,k$. Hence the second group
contains all $\s_{a_{ij}}$'s,
and one can prove inductively that the
first group
contains $\s_{a_{13}},\cdots,\s_{a_{1n}}$, hence the
second group.

\qed

\begin{lemma}
\label{quadr}
If
$a_{ij},a_{jk},a_{kl}$ and $a_{il}$ are the sides of a
topological
quadrangle taken in positive order and $a_{ik},a_{jl}$ are its
diagonals
so that $a_{ij},a_{jk},a_{ik}$ and $a_{ij}, a_{jl},
a_{il}$ are sides
of topological triangles taken in positive
order,
then
$$
\s_{a_{ij}}\s_{a_{jl}}\s_{a_{ik}} \quad =
\quad
\s_{a_{jl}}\s_{a_{ik}}\s_{a_{kl}}.
$$
\end{lemma}

\proof
This
is easily obtained from the identity
$$
\begin{array}{crcl}
&
\s_1(\s_3\inv\s_2\s_3)(\s_2\inv\s_1\s_2) & = &

(\s_3\inv\s_2\s_3)(\s_2\inv\s_1\s_2)\s_3\\
\iff &
\s_3\inv\s_1\s_2\s_3\s_1\s_2\s_1\inv & =
&
\s_3\inv\s_2\s_3\s_1\s_2\s_1\inv\s_3\\
\iff &
\s_1\s_2\s_1\s_3\s_2\s_1\inv & =
&
\s_2\s_1\s_3\s_2\s_3\s_1\inv\\
\iff & \s_2\s_1\s_2\s_3\s_2 & = &
\s_2\s_1\s_2\s_3\s_2.
\end{array}
$$
in the braid group $\br_4$,
which is homomorphically mapped
onto the group generated by
$a_{ij},a_{jk},a_{kl}$ and $a_{il}$
with
$$
\s_1\mapsto
\s_{a_{ij}},
\s_2\mapsto \s_{a_{jk}},
\s_3\mapsto
\s_{a_{kl}},
\s_2\inv\s_1\s_2\mapsto
\s_{a_{ik}},
\s_3\inv\s_2\s_3\mapsto
\s_{a_{jl}},
$$
\qed

\begin{lemma}
\label{pq-orbit}
The full-twists
of type $p$, resp.\ of type $q$ are conjugate via
twists
of type $a$ and $c$, resp.\ $b$ and $d$.
\end{lemma}

\proof
Given
$p_i,p_{i+1}$ we may apply lemma \ref{quadr},
since
$B'_i,B''_i,B'_{i+1},B''_{i+1}$ are the vertices of a
topological
quadrangle with
alternate sides $p_i,p_{i+1}$ and
$a_{i,i+1},c_{i,i+1}$.

The same argument works for $q_i,q_{i+1}$
with $b_{i,i+1},d_{i,i+1}$.

\qed

A similar argument works for the following:

\begin{lemma}
\label{su-orbit}
The
twists of type $s$ (resp.\ $u$, $u'$ or $u''$)
are conjugate via
   twists of type
$a,b,c,d$.
\end{lemma}

%%%
\subsection{Deducing the
theorem from the propositions}

\proof[Proof of theorem
\ref{bmf}]
The proof consists of three parts.

We  start  showing both inclusions
between
the braid monodromy
group and the group $H$ given in the claim;
   then we perform the weighted
count of full-twists and
their
inverses.

For the first two claims we use the braid monodromy
factorization appearing in Theorem \ref{BMF}.

We show preliminarly that
all generators of $H$ belong to the braid
monodromy group of such
a factorization.

All full-twists of type $p_1$ and $q_1$ belong to the
braid monodromy group,  by Theorem \ref{BMF}
and by our numerical assumptions.

Then the elementary factorizations $\beta_{f,i}$ and  $\beta_{g,j}$
show that  also
$\s_{a_{1i}},\s_{c_{1i}}$, resp.\
$\s_{b_{1\imath}},\s_{d_{1\imath}}$ belong to the
braid monodromy group.

In view of
lemma \ref{braid-gen} therefore all half-twists in $G$ of
type
$a,b,c,d$ belong to the braid monodromy group.

The elementary factorization $\beta_{fg,1}$ guarantees that also the
elements $\s_s,\s_{u'},\s_{u''}$ belong to the
braid monodromy group.

Concerning the other  full-twists of type $p$ and $q$, they
are in the
braid monodromy group by virtue of Lemma \ref{pq-orbit}.

For the reverse
inclusion we have to
show that $H$ contains all factors of the factorization
appearing in Theorem \ref{BMF}.

For the half twists of type $a,b,c,d$
this is immediate from lemma \ref{braid-gen}.
For the full twists of type
$p,q$ it then follows from lemma \ref{pq-orbit}.
Therefore also the half twists
of type $s$ and the cubes of half twists
of type $u,u',u''$
are in $H$
   by lemma \ref{su-orbit}.

Last we have to
perform the weighted count of full-twists
of respective types $p,q$.

It is equal
to the weighted count of intersections of $f$ with $\Phi$, resp.\
$g$
with $\Psi$, which in turn is simply the algebraic intersection
number of each pair of
divisors.  Hence
\begin{eqnarray*}
\#_p \quad =
\quad (2a,2b).(2a-c,2b-d) & = & 4ab-2ad+4ba-2bc\\
\#_q \quad = \quad
(2c,2d).(2c-a,2d-b) & = &
4cd-2cb+4dc-2da
\end{eqnarray*}
\qed

%%%%%%%%%%%%%%%%%%%%%%%%%%%%%%%
%%%%%%%%%%%%%%%%%%%%%%%%%%%%%%%%%
\section{proof
of the conjugacy theorem}

In this last section we are going to prove
that the conjugacy classes
of $\sigma_p^2$ and $\sigma_q^2$ are
different in the stabilised braid
monodromy group $\hat H$.

We rely
on the unicity of the 'roots' $\sigma_p$ and $\sigma_q$ and on
the
following proposition:

\begin{prop}
The elements $\sigma_p$ and
$\sigma_q$ are not $\hat H$-conjugate
in the braid group $G$, ie.\
there is no $\hat h\in\hat H$
such that $\hat h\sigma_p = \sigma_q
\hat h$.
\end{prop}

\Proof
The idea is to show that there is a
homomorphism
$\langle \hat H,\sigma_p,\sigma_q\rangle\to
\Sigma$,
such that the images of $\sigma_p,\sigma_q$ are not
conjugate under
the image $\hat\Ga$ of $\hat H$.

$\hat H$ consists
of liftable braids for the simple $4:1$ covering
$X_0\to \PP^1$
branched at $4(b+d)$ points associated to the
covering monodromy
homomorphism
$\pi_1(\PP^1-pts)\to \gS_4$,
but $\sigma_p$ and
$\sigma_q$ do \emph{not} lift, so $\Sigma$ can not be
taken to be
related with this covering.

Consider instead the composed
map
$\pi_1(\PP^1-pts)\to \gS_4\to \gS_3$,
and the associated $3:1$
covering $Y_0\to \PP^1$, then
the elements of $\hat H$ and
$\sigma_p,\sigma_q$ lift.

In particular we can take
$\Sigma={Sp}(H_1(Y_0,\ZZ/2))$.
Since half twists map to
transvections, full twist
map to the identity in $\Sigma$, hence the
image
$\hat\Ga$ of $\hat H$ is
the same as the image $\Ga$ of $H$.
It
thus suffices to show,
that the images of $\sigma_p$ and
$\sigma_q$
are not conjugate under the image $\Ga$ of $H$.

Since $\sigma_p$ and
$\sigma_q$ map to transvections on vectors
$v_p,v_q\in
H_1(Y_0,\ZZ/2)$ the claim follows from the following claim, which is
going to be proven in
Proposition \ref{transvections}:

CLAIM : $v_p$ and $v_q$
belong to different orbits for
the $\Ga$-action on
$H_1(Y_0,\ZZ/2)$.

\qed

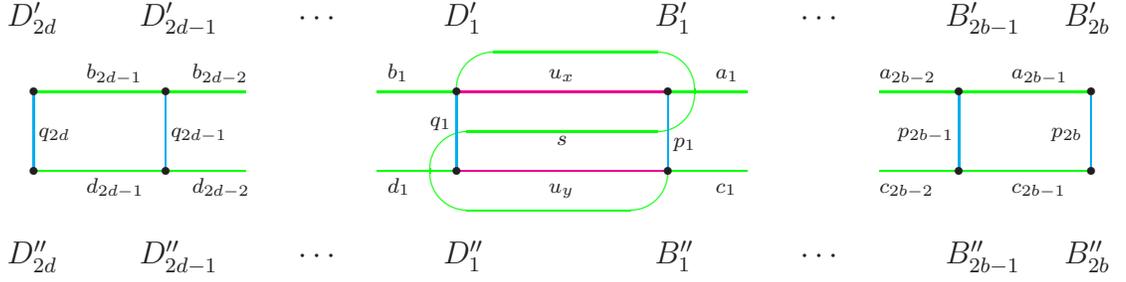
\begin{figure}
\begin{picture}(420,120)(-210,-60)

{\color{green}
\put(-200,-15){\line(1,0){80}}
\put(-70,-15){\line(1,0){30}}
\put(200,-15){\line(-1,0){80}}
\put(70,-15){\line(-1,0){30}}
\put(-200,15){\line(1,0){80}}
\put(-70,15){\line(1,0){30}}
\put(200,15){\line(-1,0){80}}
\put(70,15){\line(-1,0){30}}

\put(5,15){\oval(90,30)[t]}
\put(0,15){\oval(100,30)[rb]}
\put(0,-15){\oval(100,30)[lt]}
\put(-5,-15){\oval(90,30)[b]}
\color{magenta}
\put(-40,15){\line(1,0){80}}
\put(-40,-15){\line(1,0){80}}
\color{cyan}
\put(-200,-15){\line(0,1){30}}
\put(-150,-15){\line(0,1){30}}
\put(-40,-15){\line(0,1){30}}
\put(40,-15){\line(0,1){30}}
\put(150,-15){\line(0,1){30}}
\put(200,-15){\line(0,1){30}}
}

\put(40,15){\circle*{3}}
\put(150,15){\circle*{3}}
\put(200,15){\circle*{3}}
\put(-40,15){\circle*{3}}
\put(-150,15){\circle*{3}}
\put(-200,15){\circle*{3}}
\put(40,-15){\circle*{3}}
\put(150,-15){\circle*{3}}
\put(200,-15){\circle*{3}}
\put(-40,-15){\circle*{3}}

\put(-150,-15){\circle*{3}}
\put(-200,-15){\circle*{3}}

\put(-210,-50){$D''_{2d}$}
\put(-160,-50){$D''_{2d-1}$}
\put(-45,-50){$D''_{1}$}
\put(35,-50){$B''_{1}$}
\put(145,-50){$B''_{2b-1}$}
\put(190,-50){$B''_{2b}$}
\put(-210,40){$D'_{2d}$}
\put(-160,40){$D'_{2d-1}$}
\put(-45,40){$D'_{1}$}
\put(35,40){$B'_{1}$}
\put(145,40){$B'_{2b-1}$}
\put(190,40){$B'_{2b}$}

\put(-100,-50){$\cdots$}
\put(-100,40){$\cdots$}
\put(90,-50){$\cdots$}
\put(90,40){$\cdots$}

\put(-198,-2){$\scriptstyle
q_{2d}$}
\put(-148,-2){$\scriptstyle
q_{2d-1}$}
\put(-50,2){$\scriptstyle
q_{1}$}
\put(42,-6){$\scriptstyle p_{1}$}
\put(127,-2){$\scriptstyle
p_{2b-1}$}
\put(185,-2){$\scriptstyle
p_{2b}$}

\put(-180,20){$\scriptstyle
b_{2d-1}$}
\put(-140,20){$\scriptstyle
b_{2d-2}$}
\put(-66,20){$\scriptstyle
b_{1}$}
\put(-5,20){$\scriptstyle u_{x}$}
\put(58,20){$\scriptstyle
a_{1}$}
\put(120,20){$\scriptstyle
a_{2b-2}$}
\put(170,20){$\scriptstyle
a_{2b-1}$}
\put(-180,-23){$\scriptstyle
d_{2d-1}$}
\put(-140,-23){$\scriptstyle
d_{2d-2}$}
\put(-66,-23){$\scriptstyle
d_{1}$}
\put(-5,-23){$\scriptstyle u_{y}$}
\put(58,-23){$\scriptstyle
c_{1}$}
\put(120,-23){$\scriptstyle
c_{2b-2}$}
\put(170,-23){$\scriptstyle
c_{2b-1}$}

\put(-2,-5){$\scriptstyle
s$}
\end{picture}

\caption{isotoped generators}
\end{figure}

Let us
first consider in more detail the corresponding geometrical
properties of
the two coverings of $F_0$.

\noindent
The branched covering $X_0\to
F_0$ is associated
with a map $\pi_1(F_0-\Delta|_0)\to \gS_4$.
We
define a further branched covering $Y_0\to F_0$ by
means of the
induced homomorphism $\pi_1(F_0-\Delta|_0)\to \gS_3$
using the
natural surjection $\gS_4\to \gS_3$ with kernel given
by the Klein
four-group.

%\subsubsection{liftable twists}

For a simple
$4$-covering, taken an arc connecting two branch points,
the minimal
power of the corresponding half twist which is liftable
is
\begin{enumerate}
\item[(1)]
one, if the two associated
transpositions of the local monodromies
are equal,
\item[(2)]
two, if
they are disjoint,
\item[(3)]
three, if they do not
commute.
\end{enumerate}
We recall from theorem \ref{bmf}, that (1)
applies for arcs of type
$a,b,c,d,s$, that (2) applies for types
$p,q$ and (3) applies for type
$u$.

Considering now the associated
simple $3$-covering we note
that each arc in the base between branch
points has a
preimage, which consists of
either
\begin{enumerate}
\item[i)]
a single arc containing both
ramification points,
\item[ii)]
a cycle containing both ramification
points and a residual arc
\end{enumerate}
depending on whether the
local monodromies
with respect to a base point on the arc
yield two
different or two equal transpositions.

In fact since every branch
point corresponds to a transposition in
monodromy, the product of two
can only be a $3$-cycle or trivial
in case of a simple branched cover
of degree $3$.

Since in case ii) the half twist on the arc lifts to
a diffeomorphism
of the $3$-covering and in case i) only the cube of
the half twist does,
we have the following
implications
\begin{quote}
case $1)$ for $X_0$ $\implies$ case $ii)$
for $Y_0$,\\
case $2)$ for $X_0$ $\implies$ case $ii)$ for $Y_0$,
(if
the cube and the square\\
\hspace*{\fill} lift, then so does the half
twist itself.)\\
case $3)$ for $X_0$ $\implies$ case $i)$ for
$Y_0$.
\end{quote}
In particular only the arcs $u_i$
belong to the
first case, while all others, of types $a,b,c,d,p,q,s$,
belong to the
second case.
\\

%\subsubsection{cycles and intersection pairing}
For
the following discussion we need cycles on $Y_0$ and
their
intersection pairings.
We denote by $\tilde v$ the cycle in the
preimage of an arc $v$ in the
second class.

\begin{lem}
For the $\Z
/ 2$-valued intersection pairing $\langle\,\,,\,\rangle$ of
cycles on
$Y_0$:
$$
\langle\tilde v_1,\tilde v_2\rangle \quad\equiv_2\quad \#
(v_1\cap v_2),
$$
if $v_1\neq v_2$ and
$v_1,v_2
\in\{a_i,b_i,c_\imath,d_\imath,p_j,q_j,s\}$.
\end{lem}

\proof
It
suffices to prove that there is a bijection of points of $\tilde
v_1\cap
\tilde v_2$ with points of $v_1\cap v_2$ given by the
restriction of the
projection map.

Each branch point has two
preimages,
a ramification point and a point, where the map is
unramified,
but it is always the ramification point which belongs
to
a cycle in the preimage of any arc from that branch point.
A pair
of cycles $\tilde v_1,\tilde v_2$ meets in the ramification point
if
$v_1,v_2$ meet in the branch point.
Hence two arcs share a branch
point if and only
if they share the corresponding ramification
point.

All ordinary points have three preimages. But a cycle in
the
preimage of an arc contains either none or two of those.
Hence
either one or two of the preimages of an intersection between
two arcs
are shared by their corresponding cycles on $Y_0$.

Observe
that the only pairs of arcs which meet in a point which is not
a
branch point are $s$ and one of $a_1,d_1,p_1,q_1$.

First we consider
the four local monodromies with respect to a base
point at the
intersection of $s$ and $p_1$ along their segments.
These monodromies
are identical if and only if the cycles $\tilde s$
and $\tilde p_1$
meet in two points.
If so, we could conclude by isotopy that the two
monodromies
associated to $u_x$ must be the same.
Since on the
contrary $u_x$ belongs to class i), we have proved
that $\#\tilde
s\cap\tilde p_1=1$.

For $s$ and $a_1$ we have the same argument and
we conclude
for $s$ and $q_1$ resp.\ $d_1$ by
symmetry.
\qed

%\subsubsection{basis}

\begin{lem}
\label{hbasis}
A
basis of $H_1(Y,\Z / 2)$ is given by
$$
\tilde a_3,...,\tilde
a_{2b-1},\tilde p_{2b},
\tilde c_{2b-1},...,\tilde c_{1},\tilde
s,
\tilde b_1,...,\tilde b_{2d-1}, \tilde q_{2d},
\tilde
d_{2d-1},...,\tilde d_{2}.
$$
\end{lem}

\proof
   From the previous
lemma we can deduce that the intersection matrix
for the given
elements is tridiagonal with diagonal entries $0$
and entries in the
secondary diagonals all equal to $1$.
Hence the intersection form on
these elements is non-degenerate
and thus the elements are linearly
independent.
\\
On the other hand by the Riemann-Hurwitz formula
their number
$$
4b+4d-4 = 2b - 3 + 1 + 2b - 1 + 1 + 2d - 1 + 1 + 2d -
2
$$
equals the dimension of $H_1(Y,\Z / 2)$
since
\begin{eqnarray*}
e(Y_0) & = & 3\cdot 2-4b-4d\\
\implies\quad
b_1(Y_0) & = & 2-e \quad= \quad
4b+4d-4
\end{eqnarray*}
\qed

\begin{lem}
\label{qexists}
There is a
quadratic form ${\mathbf q}$ on $H_1$ such
that
\begin{enumerate}
\item
${\mathbf q}(w)=1$ for all elements of
the basis given in
lemma \ref{hbasis},
\item
the intersection pairing
coincides with the induced symmetric form:
$$
\langle \tilde
v_1,\tilde v_2 \rangle \quad = \quad
\langle \tilde v_1, \tilde v_2
\rangle_{\mathbf q} :=
{\mathbf q}(\tilde v_1+\tilde v_2) + {\mathbf
q}(\tilde v_1) + {\mathbf q}(\tilde
v_2).
$$
\end{enumerate}
\end{lem}

\proof
Obviously the quadratic
form associated with the matrix $Q$
meets both requirements, if the
diagonal element are
given by $Q_{ii}=1$ and
the off-diagonal
element are given by $Q_{ij}=0$ if $i>j$ and
$Q_{ij}=\langle \tilde
v_i,\tilde v_j\rangle$ if $i<j$.
\qed

We consider the subgroup
$\tilde H$ of the braid group generated
by the braid monodromy group
and by the half twists on
the arcs $p_j, q_j$.

\begin{lem}
All
elements of $\tilde H$ considered as isotopy classes
of
homeomorphisms of $(F_0,\Delta|_0)$ can be lifted to isotopy
classes
of homeomorphisms of $Y$.
\end{lem}

\proof
An element is
liftable if and only if its induced map on the natural
homomorphism
is the identity.
Since elements of the braid monodromy group
stabilise the
natural homomorphism to $\gS_4$, they stabilise also
its
composition with $\gS_4\to \gS_3$.

Half-twists on arcs $p_i$ and
$q_i$ do not stabilise the
natural homomorphism to $\gS_4$, but they
exchange
disjoint transpositions.
However such pairs are identified
under $\gS_4\to \gS_3$,
hence also the additional elements
lift.
\qed

We consider next the action of $\Ga$ via the mapping
class group
of $Y$ on the first homology of $Y$ with $\Z / 2$
coefficients
preserving the natural symplectic form.

Note that the
preimage in $Y$ of some arc in $F_0$
contains a simple closed
curve
if and only if the half-twist on that arc lifts to
$Y$.

\begin{lem}
The image of $\tilde H$ in $Sp\;
H_{\!\;\!1}\!\;\!(Y;\Z / 2)$ is generated
by the symplectic
transvection on all the classes of $H_1$
which are represented by
simple closed curves
mapping to the arcs
$a,b,c,d,q,p,s$.
\end{lem}

\proof
Since $Y\to F$ is a simple triple
cover,
any arc either lifts to an arc in $Y$ or the union of an arc
with
a simple closed curve.
In the first case the third power of the
corresponding half twist
lifts, but its lift is the half twist on the
preimage and therefore
isotopic to the identity.

In the second case
the corresponding half twist lifts to
the product of the half twist
on the arc and a Dehn twist on the circle.
The former is isotopic to
the identity but the latter acts by a
symplectic transvection by the
class of the circle.

Hence it suffices to take all transvection
associate to circles in
the lifts of the arcs needed to generated
$\Ga$.
\qed

\begin{lem}
If ${\mathbf q}$ is a quadratic form on a
$\Z / 2$-vector space
such that $\langle\, ,\,\rangle_{\mathbf q}$
is
the corresponding symmetric bilinear form,
then ${\mathbf q}(w)=1$
implies
${\mathbf q}(T_wz) \equiv_2 {\mathbf
q}(z)$.
\end{lem}

\proof
By assumption either $\langle
w,z\rangle_{\mathbf q}$ or
$\langle w,z\rangle_{\mathbf q}+{\mathbf
q}(w)$ is even, hence
\begin{eqnarray*}
{\mathbf q}(T_wz) & \equiv_2
&
{\mathbf q}\big(z+\langle w,z\rangle_{\mathbf q} w\big)\\
&
\equiv_2 &
{\mathbf q}(z)+{\mathbf q}\big(\langle w,z\rangle_{\mathbf
q} w\big)
+ \big(z,\langle w,z\rangle_{\mathbf q} w\big)\\
&
\equiv_2 &
{\mathbf q}(z) + \langle w,z\rangle_{\mathbf q}
\big(
{\mathbf q}(w) +\langle w,z\rangle_{\mathbf q}\big)\\
& \equiv_2
&
{\mathbf q}(z).
\end{eqnarray*}
\qed

\begin{lem}
If ${\mathbf
q}(w)\equiv_20$ and
$\langle w,z\rangle_{\mathbf q}\equiv_21$,
then
${\mathbf q}(T_wz) \equiv_2 {\mathbf q}(z)+1$.
\end{lem}

\proof
By
the same computation as above but with
$\langle w,z\rangle_{\mathbf
q}$
and $\langle w,z\rangle_{\mathbf q}+{\mathbf q}(w)$ both
odd.
\qed

\begin{prop}
The transvections associated to $q_{2b}$ and
$p_{2d}$
are not conjugate under the group $\Ga$ generated by the
transvections on the elements
$a_i,b_i,c_i,d_i,s$.
\end{prop}

\proof
We work with the basis of
$H_1$ given in lemma \ref{hbasis}
$$
\tilde a_3,\tilde a_4,...,\tilde
a_{2b-1},\tilde p_{2b},\tilde c_{2b-1},...,
\tilde c_2,\tilde
c_1,\tilde s,\tilde b_1,\tilde b_2,....,\tilde b_{2d-1},
\tilde
q_{2d},\tilde d_{2d-1},...,\tilde d_3,\tilde d_2.
$$
and the
quadratic form ${\mathbf q}:H_1\to\Z / 2$ of
lemma \ref{qexists},
which is non-trivial on all basis elements.

By straightforward
computation by looking at the intersection number,
we get for the
remaining elements
$a_1,a_2,d_1$:
\begin{eqnarray*}
\tilde a_1 & = &
\tilde a_3+\tilde a_5+\cdots+\tilde a_{2b-1}
+\tilde
c_{2b-1}+\cdots+\tilde c_3+\tilde c_1\\
\tilde a_2 & = &
\phantom{+1}
\tilde a_4+\tilde a_6+\cdots+\tilde a_{2b-2}+\tilde
p_{2b}
+\tilde c_{2b-2}+\cdots+\tilde c_4+\tilde c_2+\tilde s\\
& &
{}+\tilde b_2+\tilde b_4+\cdots+\tilde b_{2d-2}+\tilde q_{2d}
+\tilde
d_{2d-2}+\cdots+\tilde d_4+\tilde d_2\\
\tilde d_1 & = & \tilde
b_1+\tilde b_3+\cdots+\tilde b_{2d-1}
+\tilde d_{2d-1}+\cdots+\tilde
d_5+\tilde d_3
\\[2mm]
{\mathbf q}(a_1) & = & 1,\qquad
{\mathbf
q}(a_2) \quad=\quad0,\qquad
{\mathbf q}(d_1)
\quad=\quad1
\end{eqnarray*}
Let us denote by $\delta_p$ the element
of the dual basis, which
evaluates non-trivial on $p_{2b}$ and let us
introduce a function
$$
\rho :\,H_1(Y,\Z / 2)\to \Z / 2,\quad
x\mapsto \delta_p(x)+{\mathbf q}(x).
$$
Next we note that this
function is invariant under the action of $\Ga$,
since it is under
the
transvections associated to any one of the elements $w$ given in
the
statement of this proposition:
\begin{eqnarray*}
\rho(T_w x)) & =
& \delta_p(T_w x)+{\mathbf q}(T_w x)\\
& = & \delta_p(x)+\langle
x,w\rangle\delta_p(w)
+{\mathbf q}(x)+
{\mathbf q}\big(\langle
x,w\rangle w\big)
+\big(x,\langle x,w\rangle w\big)\\
& =
&
\delta_p(x) +{\mathbf q}(x) + \langle
x,w\rangle\big
(\delta_p(w)+{\mathbf q}(w)+\langle
x,w\rangle\big)\\
& = & \rho(x)
\end{eqnarray*}
The product on the
right in the second but last step has factors
of either parity, since
$\delta_p(w)+{\mathbf q}(w)=1$ due to
$$
\delta_p(w)\quad=\quad
\left\{
\begin{matrix}1 & \text{for }w=a_2\\ 0 &
\text{else}\hfill
\end{matrix}\right.
\qquad\qquad{\mathbf
q}(w)\quad=\quad \left\{
\begin{matrix}0 & \text{for }w=a_2\\ 1 &
\text{else}\hfill
\end{matrix}\right.
$$
Since $\rho(q)={\mathbf
q}(q)+\delta_p(q)=1+0=1$ and
$\rho(p)={\mathbf
q}(p)+\delta_p(p)=1+1\equiv_2 0$
the two elements $q,p$ are not
contained in the same $\Ga$-orbit.
\qed

\begin{prop}\label{transvections}
No full twist on any arc $p_i$ is conjugate in the braid monodromy
group to a full twist on any arc $q_j$.
\end{prop}

\proof
Suppose to the contrary that there is a pair $p_i,q_j$ such that the
corresponding full twist are conjugate under some $\beta$ from
the braid monodromy group.

Then $\beta$ maps $p_i$ to $q_j$ or vice versa.
In particular also the half-twists on $p_i$ and $q_j$ are
conjugate under $\beta$ in the braid monodromy group.

But then the transvections associated to $p$ and $q$
are conjugate under the symplectic transformation associated
to $\beta$.
That is in contradiction to what we proved earlier,
so our claim must be true.
\qed

{\bf Acknowledgements.}
The present research was performed in the realm of the Forschergruppe 790
`Classification of algebraic surfaces and compact complex manifolds'
of the DFG (DeutscheForschungsGemeinschaft).

The results of the paper were announced in a workshop which took place in
october 2008 during the special trimester "Groups in algebraic Geometry"
at the Centro De Giorgi of the Scuola Normale di Pisa.

\bigskip

%%%%%%%%%%%%%%%%%%%%%%%%%%%%%%%%%%%

\noindent
{\bf Author's addresses:}

\bigskip

\noindent
Prof. Fabrizio Catanese\\
Lehrstuhl Mathematik VIII\\
Universit\"at Bayreuth, NWII\\
         D-95440 Bayreuth, Germany

e-mail: Fabrizio.Catanese@uni-bayreuth.de

\noindent
Prof. Michael L\"onne\\
Lehrstuhl Mathematik VIII\\
Universit\"at Bayreuth, NWII\\
         D-95440 Bayreuth, Germany

e-mail: Michael.Loenne@uni-bayreuth.de

\noindent
Prof. Dov Bronislaw Wajnryb \\
Department of Mathematics and Applied Physics\\
Technical University of Rzeszow\\
        Rzeszow, Poland

e-mail : dwajnryb@prz.edu.pl

\end{document}